\documentclass{aims} 

\usepackage{amsmath}
\usepackage{amssymb}
\usepackage[utf8]{inputenc}

\DeclareUnicodeCharacter {0010}{C}
\usepackage{epsfig} 
\usepackage{xcolor}
\usepackage{bm}
\usepackage{epstopdf} 
\usepackage[colorlinks=true]{hyperref}
\hypersetup{urlcolor=blue, citecolor=red}
\allowdisplaybreaks

\def\currentvolume{X}
 \def\currentissue{X}
  \def\currentyear{200X}
   \def\currentmonth{XX}
    \def\ppages{X--XX}
     \def\DOI{10.3934/xx.xxxxxxx}

\newtheorem{theorem}{Theorem}[section]
\newtheorem{corollary}{Corollary}

\newtheorem{lemma}[theorem]{Lemma}

\theoremstyle{definition}
\newtheorem{definition}[theorem]{Definition}
\newtheorem{remark}{Remark}

\title[Persistent homology with k-nearest-neighbor filtrations]
{Persistent homology with k-nearest-neighbor filtrations reveals topological convergence of PageRank} 

\author[MINH QUANG LE and DANE TAYLOR]{}


 \email{minhquan@buffalo.edu}
 \email{danet@buffalo.edu}

\thanks{The first author is supported by NSF grant xx-xxxx}

\begin{document}
\maketitle
\centerline{\scshape Minh Quang Le$^*$ and Dane Taylor$^*$}
\medskip
{\footnotesize
 \centerline{Department of Mathematics}
   \centerline{  University at Buffalo, State University of New York }
   \centerline{USA}
} 

\medskip


\bigskip

 \centerline{(Communicated by the associate editor name)}


\begin{abstract}
Graph-based representations of point-cloud data are widely used in data science and machine learning, including $\epsilon$-graphs that contain edges between pairs of data points that are nearer than $\epsilon$ and kNN-graphs that connect each point to its $k$ nearest neighbors.
Recently, topological data analysis has emerged as a family of mathematical and computational techniques to investigate topological features of data using simplicial complexes. These are a higher-order generalization of graphs and many techniques such as Vietoris-Rips (VR) filtrations are  also parameterized by a distance $\epsilon$.
Here, we develop kNN complexes as a generalization of kNN graphs, leading to kNN-based  persistent homology techniques for which we develop stability and convergence results.
We apply this technique to characterize the convergence properties  PageRank, highlighting how the perspective of discrete topology complements traditional geometrical-based analyses of convergence. Specifically, we show that   convergence of relative positions (i.e., ranks) is captured by kNN persistent homology, whereas  persistent homology with VR filtrations   coincides with vector-norm convergence.
Beyond PageRank, kNN-based persistent homology is expected to be useful to other data-science applications in which the relative positioning of data points is more important than their precise locations.

\end{abstract}


\section{Introduction}

Topological data analysis (TDA) is a rapidly growing field of applied mathematics in which techniques from computational and applied topology are applied to extract structural  information about the ``shape'' of data. TDA has been applied to numerous contexts ranging from 
visualization and dimensionality reduction \cite{carlsson2009topology,kusano2016persistence} and  time series analyses \cite{perea2015sliding} to applications  in
cosmology \cite{sousbie2011persistent,weygaert2011alpha}, 
physical processes over networks
\cite{le2021persistent,taylor2015topological,kramar2013persistence,kondic2012topology}, neuroscience \cite{petri2014homological,giusti2015clique,chung2009persistence,nielson2017uncovering}, and systems biology \cite{rucco2014using,liang1998analytical,kasson2007persistent,ichinomiya2020protein}. 
One of the main tools is the study of  persistent homology, which can effectively reveal multiscale topological properties of data. This approach relies on examining how the homology of a topological space evolves as one applies a filtration, which is one of the basic notions in topology. There are many types of filtrations \cite{aktas2019persistence,rieck2017clique,huang2016persistent,chowdhury2018functorial,gasparovic2018relationship}, however we will focus on the widely used Vietoris-Rips  (VR) filtration \cite{edelsbrunner2008persistent}.
In general, different filtrations reveal complementary insights, and it is important to develop additional filtrations that cater to different applications.

\begin{figure}[t] 
\begin{center}
    \includegraphics[width=\linewidth]{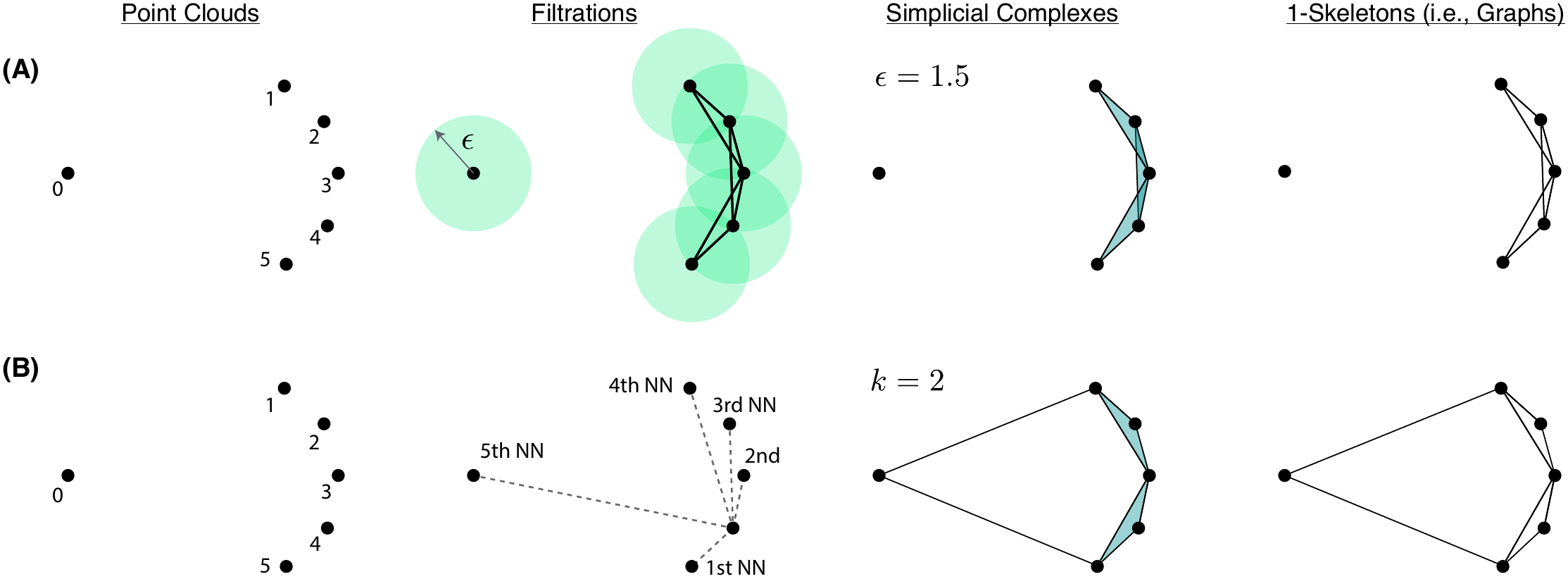}
    \caption{
    Visualizations of simplicial complexes and graphs resulting from (A) a Vietoris-Rips (VR) filtration and (B) our proposed k-nearest-neighbor (kNN) filtration. As shown in the second column, VR filtrations are parameterized by the radius $\epsilon$ of $\epsilon$-balls that are centered at the points, whereas kNN filtrations are parameterized by the number $k$ of nearest neighbors. (Dotted lines depict the nearest-neighbor orderings for node $i=4$.) The third column depicts simplicial complexes that are obtained at   some $\epsilon$ and $k$. The fourth columns shows their 1-skeletons, which are graphs in which $k$-simplices of $k>1$ are discarded.
    }
    \label{fig:fig1}
    \end{center}

\end{figure}

In the prototypical setting, one aims to construct and study empirical topological features for  a set of data points, or \emph{point cloud}.  The approach involves using points to construct a filtered topological space, which is often represented by a simplicial complex in which the data points are 0-simplices, and higher-dimensional $k$-simplices are constructed via some set of rules. Often, $k$-simplices are added according   to the pairwise distances between 0-simplices. 
As an example, in Fig.~~\ref{fig:fig1}(A), we visualize one of the most  commonly studied point-cloud filtrations, the Vietoris-Rips (VR) filtration. As shown, VR filtrations are constructed by considering simplicial complexes in which k-simplices are added between 0-simplices that are less that $\epsilon>0$ distance apart.

In Fig.~~\ref{fig:fig1}(B), we illustrate a different filtration that we  develop in this paper: \emph{kNN filtrations}, whereby $k$-simplices are created according to $k$-nearest-neighbor sets.  We note that while $\epsilon$ and kNN graphs are both very prevalent in the data-science and machine-learning literatures, TDA methods  largely focus on  $\epsilon$-based simplicial complexes and filtrations, thereby limiting their potential utility to new applications.
We develop kNN-based simplicial complexes as a generalization of kNN graphs,  allowing us  to develop kNN-based   TDA methods including  kNN persistent homology. 
%
We formulate and study persistent homology under kNN filtrations,  constructing persistence diagrams and studying their robustness properties. We also define and study several local version of kNN complexes, filtrations, and persistent homology. By constructing filtrations using discrete sets, our approach relies on discrete topology, thereby contrasting   approaches that are tied to continuous topological spaces \cite{kim2019homotopy}. 
In addition, we further analyze and explore  homological features for converging sequences of point sets, exploring how the convergence of  persistent diagrams for kNN filtrations   contrasts that for VR filtrations. 

\begin{figure}[t]
\begin{center}
\includegraphics[width=.42\linewidth]{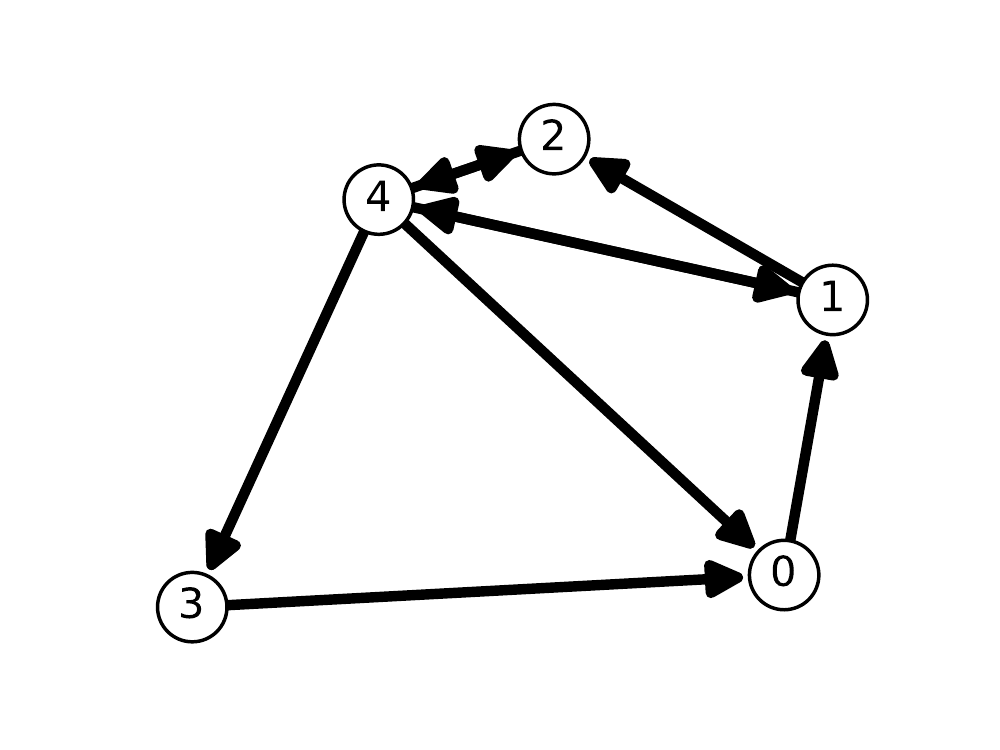}
\includegraphics[width=.5\linewidth]{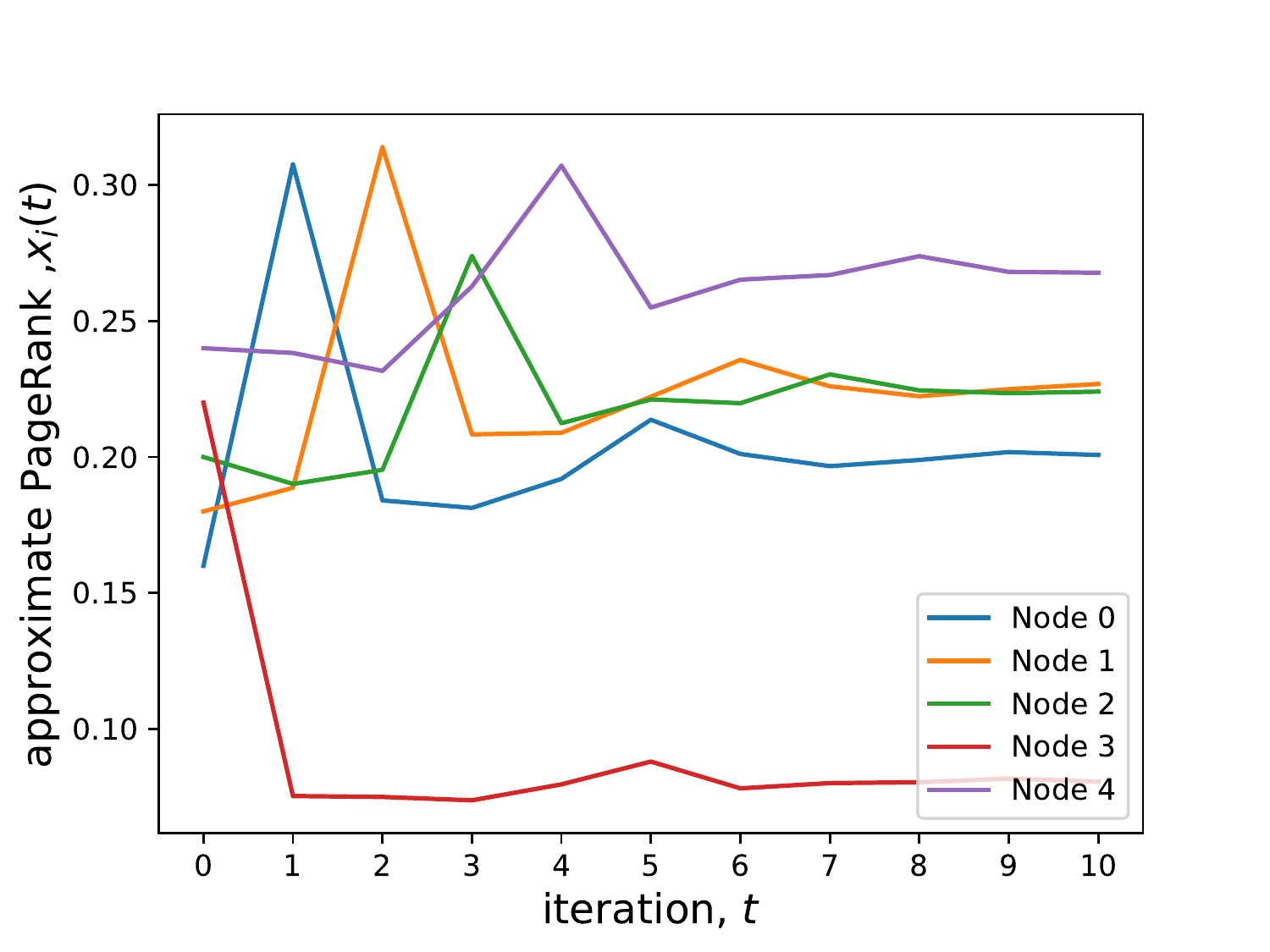}
\caption{ 
(left)~An example graph.
(right)~Convergence  of  approximate  PageRank values $x_i(t)\to\pi_i$ with $t$ iterations. 
}
\label{Fig:graphPR}
\end{center}
\end{figure}

We   apply kNN persistent homology to study the ranking of nodes in graphs, and in particular, we study the convergence of rankings given by approximate PageRank values that are obtained by the power iteration method.
Our work establishes a new connection between TDA and PageRank, and provides a topological approach for determining how many iterations are required for the node rankings to convergence.
%
We visualize this application and motivation in Fig.~\ref{Fig:graphPR}. In this application, approximate PageRank  values $x_i(t)\to\pi_i$ asymptotically converge to their final values and the  normed error $||{\bf x}(t) -\bm{\pi}||$ exponentially decays. However, as highlighted in Table \ref{table:student}, the orderings (i.e., node ranks)  according to ${ x}_i(t)$ values converge   after only $t=9$ iterations, even though $||{\bf x}(t) -\bm{\pi}||>0$ for all $t$. 
We apply kNN-based persistent homology to study of convergence for these relative orderings, which is a property that is crucial to PageRank and which is not  revealed through VR filtrations.
Although we focus here on PageRank, we expect kNN-based TDA  to be widely useful to other applications in which the relative positioning of points---as opposed to the precise locations---is a property of main interest.



\begin{table}[t]
\centering
{\tiny{
\begin{tabular}{ | l | r | r | r | r | r | r | r | r | r | r | r |}
    \hline
    Nodes & $R_i(0)$ & $R_i(1)$ & $R_i(2)$ & $R_i(3)$ & $R_i(4)$ & $R_i(5)$ &$R_i(6)$& $R_i(7)$&$R_i(8)$&$R_i(9)$&   ${R}_i(\infty)$ \\ \hline \hline
     $i=0$ & 5 & 1 & 4 & 4 & 4 & 4 & 4 & 4 & 4 &4 &4 \\ \hline
     $i=1$ & 4 & 4 & 1 & 3 & 3 & 2 & 2 & 3 & 3 &2 &2   \\ \hline
     $i=2$ & 3 & 3 & 3 & 1 & 2 & 3 & 3 & 2 & 2 &3 &3 \\ \hline
     $i=3$ & 2 & 5 & 5 & 5 & 5 & 5 & 5 & 5 & 5 &5 &5 \\ \hline
     $i=4$ & 1 & 2 & 2 & 2 & 1 & 1 & 1 & 1 & 1 &1 &1 \\ \hline
    \end{tabular}
    \caption{The relative orderings (i.e., node ranks $R_i(t)$)  converge at $t\ge 9$. 
    For each $t$, we indicate the top-ranked node by $R_i(t)=1$. 
}}}
    \label{table:student}
\end{table}


This paper is organized as follows.
We provide background information in Sec.~\ref{sec:prior}  and our main findings in Sec.~\ref{sec:main}.
In Sec.~\ref{sec:PagerankandRD}, we apply this approach to study convergence for the PageRank  algorithm. We provide a discussion in Sec.~\ref{sec:discusion}.

\section{Background information}\label{sec:prior}

In Sec.~\ref{sec:simplical}, we discuss simplicial complexes that are derived from point clouds.
In Sec.~\ref{sec:holomology}, we define homology and persistent homology for simplicial complexes.
In Sec.~\ref{sec:stability}, we discuss the notions of stability and  convergence   for  persistence diagrams.

\subsection{Simplicial complexes derived from point clouds} \label{sec:simplical}

\subsubsection{Simplicial complexes}

We will define simplicial complexes in a geometric way by considering a  set of points, or ``point cloud'' \cite{schaub2018random}.
Consider a finite set  $ \mathcal{Y}=\{y_i \} \in\mathbb{R}^p$ of $N$ points (which  we enumerate $i\in\mathcal{V}=\{1,\dots,N\}$) in a $p$-dimensional space that we equip with the Euclidean metric. 

\begin{definition}[Euclidean Metric Space]
\label{def:Euclidean}
Let  $\mathbb{R}^p$ be the set of all ordered $p$-tuples, or vectors, over the real numbers and $d:\mathbb{R}^p\times \mathbb{R}^p$ be the Euclidean metric 
$$d(x,y)=||x-y||_2=\sqrt{\sum_{i=1}^p(x_i-y_i)^2},$$
where $x=(x_1,..,x_p), y=(y_1,..,y_p) \in \mathbb{R}^p$.
Then the Euclidean metric space is given by $(\mathbb{R}^p,d)$.
\end{definition}
To facilitate later discussion, it is also helpful to define the following  distance-related concepts.

\begin{definition}[Euclidean Ball, or $\epsilon$-ball]
\label{def:EucBal}
A $p$-dimensional Euclidean ball $\mathcal{B}_\epsilon(x)$ centered at $x$ with radius $\epsilon$ is defined by
$$\mathcal{B}_\epsilon(x)=\{y\in \mathbb{R}^p:d(x,y) \le \epsilon\}$$
\end{definition}

\begin{definition}[Pairwise-Distance Map] 
\label{def:NNOF}
Let $\mathcal{Y} =  \{y_i  \}_{i\in\mathcal{V}} \subset \mathbb{R}^p$ be a finite  set of $N$ points in a Euclidean metric space $(\mathbb{R}^p,d)$. 
Then we define the ``pairwise-distance map'' 
$f:\mathcal{Y} \to \mathbb{R}_+^{N\times N}$ 
to encode the distances between all pairs of points so that
$f_{ij}(\mathcal{Y}) = d(y_i,y_j) $. 
We similarly define a map $f_i:\mathcal{Y}\to \mathbb{R}_+^N$ for each row $i$ of matrix $[f_{ij}]$ by $f_i:\mathcal{Y}\to \mathbb{R}_+^N$ so that $[f_i(\mathcal{Y})]_j = f_{ij}(\mathcal{Y})$.
\end{definition}

Given a point-cloud, we define a simplicial complex   as a collection of $k$-simplices involving a set $ \mathcal{V}=\{1,\dots,N\}$ of vertices with locations $\mathcal{Y}=\{y_i \}$.
A $k$-dimensional simplex---or simply \textbf{k-simplex}---$\mathcal{S}^k$ is defined by a subset of $\mathcal{V}$ having cardinality $k+1$. For example, each 0-simplex is a vertex, each 1-simplex is an edge, 2-simplices are informally ``triangles'' that must be defined using $k+1=3$ vertices. Because the vertices are spatially embedded in $\mathbb{R}^p$, each $k$-simplex is a $k$-dimensional geometrical object. It is then natural to consider their $(k-1)$-dimensional faces, which can be defined as follows.  A face of a $k$-simplex $S^k$ is a subset of $S^k$ with cardinality $k$, i.e , with one of the elements of $S^k$ omitted. If $S^{k-1}_f$ is a face of simplex $S^k$, then $S^k$ is called $\textit{coface}$ of $S^{k-1}_f$. For example, the 0-simplices at the start and end of a 1-simplex are its faces, and the set of 1-simplices that are adjacent to a 0-simplex are its cofaces. 

A   simplicial complex $X$ is a collection of $k$-simplices with the property that if $S^k \in X$, then all the faces of $S^k$ are also members of $X$.
The notions of face and coface can also be extended to \textit{abstract simplicial complexes} that   lack spatial coordinates, and one can also define connections between $k$-simplices of the same dimension  by considering if their faces and cofaces overlap. Two $k$-simplices $S^k_i$  and $S^k_j$  are said to be lower adjacent if they have a common face, and they are upper adjacent if they are both faces of a common $(k + 1)$-simplex. For any $S^k \subset X$  we define its degree, denote by $deg(S^k)$, to be the number of cofaces of $S^k$. We use $X^k$ to denote the subset of $k$-simplices in $X$.

Note that a graph, while typically defined via two sets (vertices, edges), can also be interpreted as a 1-dimensional abstract simplicial complex, since it contains $k$-simplices of dimension $k \leq 1$. Moreover, for any simplicial complex, one can obtain an associated graph that is its 1-skeleton, whereby one discards any $k$-simplices of dimension $k>1$. More generally, a $\kappa$-skeleton of a simplicial complex $X$ can be obtained by discarding all $k$-simplices of dimension $k>\kappa$.
Thus, a simplicial complex can   be understood as a generalization of a graph that allows for higher-order relationships between vertices. To emphasize this connection, we will make no distinction between 1-simplices in a simplical complex $X$ and the edges of a graph, and we will refer to them interchangeably.

\subsubsection{Two simplicial complexes based on   $\epsilon$}

There are many ways in which one can construct a set of simplices involving the vertices $\mathcal{V}$ associated with a given a set of points $ \mathcal{Y}=\{y_i \} \in\mathbb{R}^p$. Most approaches stem from considering $\epsilon$-balls centered at the points $\mathcal{Y}$.
Here, we present two closely related simplicial complexes that are parameterized by a distance threshold $\epsilon$ and are often motivated by the  assumption that the point cloud lies on a low-dimensional manifold.

\begin{definition}[\v Cech Complex \cite{bubenik2015statistical}]\label{def: Cech complex}
Given a collection of points $  \mathcal{Y}=\{y_i \}_{i=1}^N \in \mathbb{R}^p$, the {\v C}ech complex, $\mathcal{C}_\epsilon$, is the abstract simplicial complex whose $k$-simplices are determined by the unordered $(k + 1)$-tuples of points $\{y_i\}_{i=1}^{k+1}$  whose closed $(\epsilon/2)$-ball neighborhoods have a point of common intersection.
\end{definition}

The $\textit{Vietoris-Rips (VR) complex}$ is   closely related and is defined as follows.

\begin{definition}[Vietoris-Rips Complex {\cite{attali2013vietoris}}]
\label{def:Vietoris-Rips complex}
Given a collection of points $\mathcal{Y}=\{y_i \}_{i=1}^N\in \mathbb{R}^p$, the Rips complex, $\mathcal{R}_\epsilon$, is the abstract simplicial complex whose $k$-simplices correspond to unordered $(k + 1)$-tuples of points $C = \{y_i\}^{k+1}_{i=1} \in \mathcal{R}_\epsilon$ whose pairwise distance satisfy  $|| y_i-y_j||_2 \le \epsilon $ for any $y_i,y_j\in C$.
\end{definition}

Notably, \v Cech complexes are sometimes preferred over VR complexes because of their relation to the a continuous topological space that is associated with   unions of $\epsilon$-balls, as established in the following Theorem.

\begin{theorem}[\v Cech theorem {\cite{kim2019homotopy}}]
The $\textbf{\v Cech theorem}$ (or, equivalently, the “Nerve theorem”) states that a {\v Cech} complex $\mathcal{C}_\epsilon$  has the homotopy type of the union of closed  $(\epsilon/2)$-balls about the point set $\mathcal{Y}=\{y_i \}$.
\end{theorem}

This result is somewhat surprising since a {\v Cech complex} $\mathcal{C}_\epsilon$ is an abstract simplicial complex, which is a discrete topological space of potentially high dimension. In contrast, the union of closed  $(\epsilon/2)$-balls is a continuous topological space (i.e., subset of $\mathbb{R}^p$),  and so it is interesting that these two spaces have the same homotopy type. While an identical relation has not been identified for VR complexes, it has been shown that they are very closely related {\v Cech complexes} through the following interleaving result.

\begin{lemma}[{Relation between \v Cech complex and Vietoris-Rips complex \cite{kim2019homotopy} }]
For any $\epsilon > 0$, there is a chain of inclusion maps for {\v Cech} complexes and VR complexes
\[
\mathcal{R}_\epsilon \subseteq \mathcal{C}_{\epsilon \sqrt{2}} \subseteq \mathcal{R}_{\epsilon \sqrt{2}}.
\]
\end{lemma}

%

From a practical perspective, applications involving TDA often focus on VR complexes because they are more computationally efficient to compute than 
{\v Cech complexes}. That is, it is easier to compute whether the distance between two points is less than $\epsilon$ versus check whether the intersection of $\epsilon$-balls is nonempty. As such, TDA methods based on VR complexes  are very popular in applied settings, and we will later conduct numerical experiments that use them as a baseline for comparison. Given the relations between VR and  {\v Cech complexes} and the unions of $\epsilon$-balls, these all may approximate the structure of a manifold, assuming all data points lie on the manifold. Such an assumption, however, is not appropriate for every data set.

Before continuing, we note that simplicial complexes can be constructed in ways that do not require a set of points. For example, one can construct abstract simplicial complexes based on a given graph.

\begin{definition}[Clique Complex \cite{zomorodian2010fast}] \label{def:Clique_complex}
The  $\textbf{clique complex}$  $Cl(G)$  of an undirected graph $G=(\mathcal{V},\mathcal{E})$  is a simplicial complex where  $\mathcal{V}$ are vertices of $G$   and each $k$-clique (i.e. a complete subgraph with $k$ vertices) in $G$ corresponds to a $(k-1)$-simplex in $Cl (G)$. 
More precisely, it is the simplicial complex
\[
Cl(G) =\mathcal{V} \cup \mathcal{E} \cup \left\{\sigma | \binom{\sigma}{2} \subseteq \mathcal{E}\right\},
\]
where $\sigma \subset \mathcal{V} $ denotes a simplex
and $\binom{\sigma}{2} \in \mathcal{V}\times \mathcal{V}$ denotes the set of pairs that can be obtained by selecting two 0-simplices from $\sigma$ (which must be an edge  in the graph).
\end{definition}

It is worth noting that the \emph{1-skeleton} (i.e., 0-simplices and 1-simplices) of any clique complex recovers its associated   graph $G=(\mathcal{V},\mathcal{E})$. That is, the construction of a graph's clique complex is an invertible transformation.

\subsection{Homology and persistent homology}
\label{sec:holomology}

After representing a data set by  a simplicial complex, one can study its topological structure via its homology. To obtain multiscale insights, it is also useful to study how that homology changes and persists as one varies a parameter, such as a distance threshold $\epsilon$.

\subsubsection{Homology}

We formulate homology by considering functions defined over $k$-simplices within a simplicial complex.
A chain complex over a field $\mathbf{F} $ is a tuple  $(\mathcal{C}; \partial)$  where $\mathcal{C}$ is a collection  $\{\mathcal{C}_{k}\}_{k\in \mathbf{N}}$ of vector spaces together with a collection $\partial$ of $\textit{F}$-linear maps $\{\partial_k: C_k \longrightarrow C_{k-1}\}_{k \in \mathbf{N}}$  such that  $\partial_{k-1} \circ \partial_k =0$ for all integers $k$. The maps $\partial_k$ are called boundary maps. The $k$-cycles of the complex are the elements that are sent to zero by the map $\partial_k$; the $k$-boundaries are the elements in the image of $\partial_{k+1}$. A map of chain complexes $f : (C; \partial) \longrightarrow (C^\prime; \partial^\prime)$ is a collection $\{f_k: C_k \longrightarrow C_{k^\prime}\}_{k \in \mathbf{N}}$  of $\mathbf{F}$-linear map such that $f_{k-1} \circ \partial_k = \partial
_k^{\prime} \circ f_k$  for all natural numbers $k$ .

The $k$-cycles form a vector space, 
and so do the $k$-boundaries; we denote these vector spaces by $Z_k$ and $B_k$, respectively. The $k$th homology of a chain complex $(C, \partial)$ over a field $\mathbf{F}$ is the quotient vector space 
$$H_k((C, \partial)) = Z_k/B_k.$$ 
The number
\[
 \beta_k (C) \equiv \textsf{dim}(H_k(C,\partial))=\textsf{dim} Z_k- \textsf{dim} B_k
\]
is called the $k$-th Betti number. For a geometric intuition, the dimension of $H_k(X)$ can be thought of as the number of `$k$-dimensional holes' of X:
\begin{itemize}
    \item  The 0th Betti number $\beta_0$ is the number of connected components.
    \item  The 1st Betti number $\beta_1$ counts the number of loops, or 1-cycles.
    \item  The 2nd Betti number $\beta_2$  counts the number of voids, or 2-cycles.
\end{itemize}
The dimension of a simplicial complex is the maximum over the dimensions of its simplices. If $X$ is a simplicial complex of dimension $d$ then $H_k(X) = 0$ for all $k \geq d.$

Any map of chain complexes $\Xi : (C, \partial) \longrightarrow (C^\prime; \partial^\prime)$ induces a linear map on homology:
\[
H(\Xi): H_k((C,\partial)) \longrightarrow H_k((C^\prime,\partial^\prime))
\]
To simplify   notation, we later use $\Xi$ in place of $H(\Xi)$.

\subsubsection{Filtrations}

To discuss how homology changes and persists, it is necessary to define a filtration for simplicial complexes.
Given a simplicial complex $X$, a $\textit{filtration}$ is a totally ordered set of subcomplexes $X_i$  of $X$, indexed by the nonnegative integers, such that if $i \leq j$  then $X_i \subset X_j$ \cite{wang2012basic} (or equivalently a sequence of simplicial complexes $ X_1,...,X_l$ such that $X_1 \subseteq \dots  \subseteq X_l = X$). The totally ordering itself is called a $\textit{filter}$. We call a simplicial complex together with a filtration a $ \textit{filtered simplicial complex}$.

There are different ways to define a filtration \cite{aktas2019persistence}, and we will focus herein on
the popular $\textit{Vietoris-Rips filtration}$.

\begin{definition}[Vietoris-Rips Filtration \cite{aktas2019persistence}]\label{VR filtration} 
Let $G=({V},{E})$  be an undirected graph 
in which each vertex $i\in\mathcal{V}$ has a position $x_i\in \mathcal{R}^p$ in a metric space with metric $d$.
Consider a weight function $W: \mathcal{V} \times \mathcal{V} \longrightarrow \mathbb{R}$ defined on the edges $\mathcal{E}$ so that $W(i,j)=d(x_i,x_j)$ encodes the distance between  $x_i$ and $x_j$. Let $Cl(G)$ be the associated clique complex for $G(\mathcal{V},\mathcal{E})$. For any $\delta \in \mathbb{R}$, the 1-skeleton $G_\delta = (V,E_\delta)$ is defined as the subgraph of $G$ where $E_\delta \subseteq E$   includes only the edges such that  $d(x_i,x_j)\le \delta$. Then, for any $\delta \in \mathbb{R}$, we define the $\textit{Vietoris-Rips filtration}$ by
\[
\{Cl(G_\delta) \longrightarrow Cl(G_{\delta^ \prime})\}_{0\leq \delta \leq \delta^\prime} .
\]
\end{definition}

In this work, we will only consider when the distance function is the Euclidean metric, although one  can in principle use other metrics. That said, we note in passing that one can also construct filtrations in a variety of ways including functional metric filtrations  \cite{aktas2019persistence}, vertex-based clique filtrations \cite{rieck2017clique}, $k$-clique filtrations  \cite{rieck2017clique}, weighted simplex filtrations  \cite{huang2016persistent}, vertex function based filtrations  \cite{aktas2019persistence}, Dowker sink and source filtration \cite{chowdhury2018functorial}, and the intrinsic \v Cech filtrations  \cite{gasparovic2018relationship}.

\subsubsection{Persistent homology}

Persistent homology involves studying how homology changes (or more precisely, identifies when it does not change) during a filtration.
Consider a simplicial complex, for all $i \leq j$ the inclusion maps $X_i \hookrightarrow {X_j}$ induce $\mathbf{F}$-linear maps $\Xi_{i,j} : H_k(X_i) \longrightarrow H_k(X_j)$ on simplicial homology. For a generator $x \in H_k(X_i)$ with $x \neq 0$, we say that $x$ dies in $H_k(X_j)$ if $j > i$ is the smallest index for which $\Xi_{i,j}(x) = 0$. 
We similarly say that $ x \in H_k(X_i)$ is born in $H_k(X_i)$ if $\Xi^{-1}_{t,i} (x) = 0$ for all $t < i$. We can represent the lifetime of $x$ by the half open interval $[i, j)$. If $\Xi_{i,j}(x) \neq 0$ for all $i < j \leq l$, then we say that $x$ lives forever and we represent its lifetime by the interval $[i;\infty)$.

\begin{figure}[t!]
    \centering
    \includegraphics[width=1\linewidth]{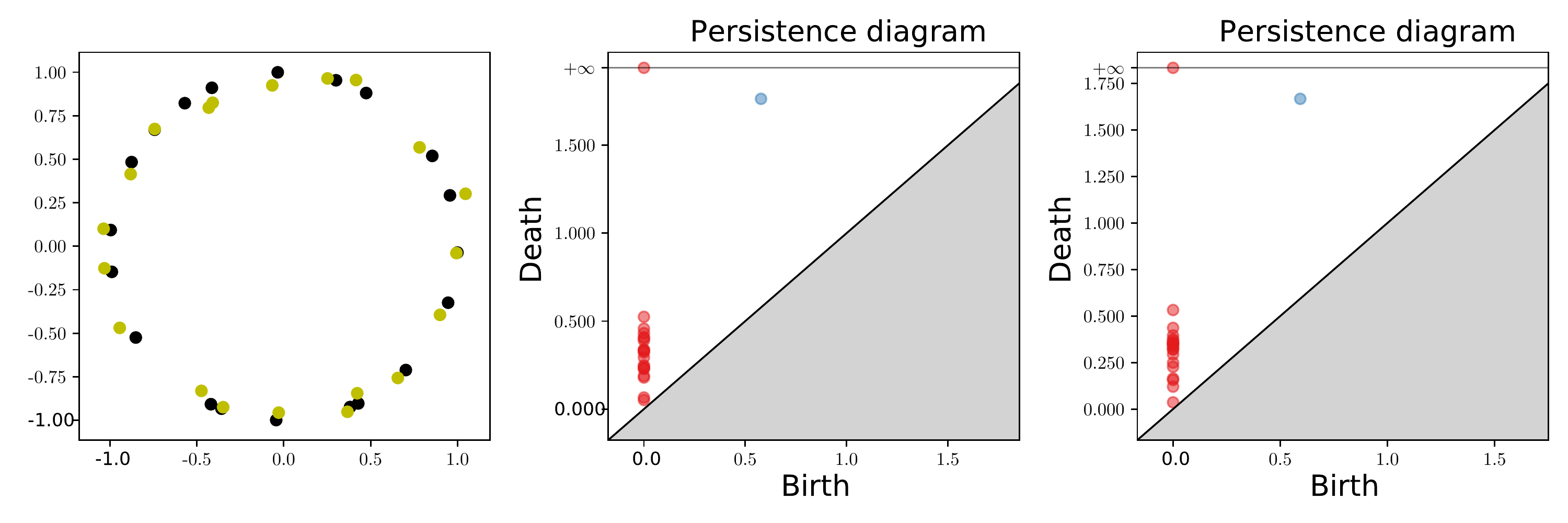}
    \vspace{-.1cm}
    \caption{Example of stability for persistence diagrams resulting from Vietoris-Rips filtrations.
    (left) Two point-cloud sets $\mathcal{Y}$ and $\mathcal{Z}$  are close with respect to the $\mathsf{L}_\infty$ norm. The center and right panels depict their associated persistence diagrams $D_{\mathcal{Y}}$ and $D_{\mathcal{Z}}$, which are also close with respect to the bottleneck distance.
    }
    \label{fig:SPD}
\end{figure}

The $k$-th persistent homology vector spaces of a filtered simplicial complex $X$ are defined as $H_k^{i,j}= \mathbf{img}(\Xi_{i,j})$, and the total $k$th persistent homology  of $X$  is defined as $\oplus _{i=1}^l H_k(X_i)$. By the Correspondence Theorem of Persistent Homology \cite{zomorodian2005computing}, for each $k$  we can assign to the total $k$-th persistent homology vector space a finite well-defined collection of half open intervals, i.e., its so called \emph{barcode}. An alternative way to represent persistent homology graphically is given by persistence diagrams, in which case each open interval $[i, j)$ is represented by the point $(i, j)$ in $\mathbb{R}^2$. See Fig.~\ref{fig:SPD} for an example persistence diagram, which we will more formally define   in the next section.

\subsection{Stability and convergence of persistence diagrams}\label{sec:stability}
Persistence diagrams and barcodes are widely used as a concise representation for the multiscale homological features of a point-cloud data set. As such, it is important to understand their robustness to perturbations, which can represent data error or noise \cite{cohen2007stability}. One can also study  for a converging sequence of point clouds whether their associated persistence diagrams also converge. In this section, we present basic results related to stability and convergence.

\subsubsection{Stability}
We first discuss stability in a general setting, which requires several definitions. Let $X$ be a simplicial complex and $g : X \longrightarrow \mathbb{R}$ be a tame function.

\begin{definition}[Tame \cite{cohen2007stability}] \label{def:tame}
A function $g : X \longrightarrow \mathbb{R} $ is tame if it has a finite number of homological critical values and the homology groups $H_k (g^{-1}(-\infty,a))$  are finite-dimensional for all $k \in \mathbb{Z}$ and $a \in \mathbb{R}$.
\end{definition}
In particular, Morse functions on compact manifolds are tame, as well as piece-wise linear functions on finite simplicial complexes and, more generally, Morse functions on compact Whitney-stratified spaces \cite{goresky1988stratified}.
 Assuming a fixed integer $k$, we define  $G_x = H_k(g^{-1}(-\infty,x])$, and for any $x<y$, we let $g_x^y:G_x \longrightarrow G_y$  be the map induced by inclusion of the sub-level set of x in that of y.

\begin{lemma}[Critical Value Lemma \cite{cohen2007stability}]\label{lemma:critical}
If some closed interval $[x, y]$ contains no homological critical value of $g$, then $g^y_x$ is an isomorphism for every integer $k$.
\end{lemma}

We now   formally define a $\textbf{persistence diagram}.$ Using the same notation as above, we write $G_x^y= {\bf img}( g_x^y) $  for the image of $G_x$  in $G_y$. By convention, we set $G_x^y=\{0\}$ whenever $x$  or $y$  is infinite. The group $G_x^y$ is called the persistent homology group.

Let $g:X \longrightarrow \mathbb{R}$  be a tame function and denote $(a_i)_{i=1,...,n}$ as its homological critical values, and let $(b_i)_{i=0..n}$ be an interleaved sequence in which  $b_{i-1} < a_i < b_i$ for all $i$. We set $b_{-1} = a_0 = -\infty$ and $b_{n+1} = a_{n+1} = +\infty$. For two integers $0 \leq i < j \leq n + 1$, we define the multiplicity of the pair $(a_i, a_j)$ by:
\[
\mu_i^j=\beta_{b_{i-1}}^{b_j}-\beta_{b_{i}}^{b_j}+\beta_{b_{i}}^{b_{j-1}}-\beta_{b_{i-1}}^{b_{j-1}}
\]
where $\beta_x^y = {\bf{dim}}( F_x^y)$ denotes the persistent Betti numbers for all $-\infty \leq x\leq y \leq +\infty$.

\begin{definition}[Persistence Diagram \cite{cohen2007stability}]\label{def:Persistence diagram} 
The persistence diagram $D(\mathcal{A}) \subset \overline{\mathbb{R}}^2$ of a filtered topological space $\mathcal{A}$  is the set of points $(a_i, a_j)$, which are counted with multiplicity $\mu_i^j$ for $0\leq i < j \leq n+1$, along with all points on the diagonal, which are counted with infinite multiplicity.
\end{definition}

Given two persistence diagrams, we can quantify their dissimilarity using various metrics. In particular, we will consider the following   metric.
 
\begin{definition}[Bottleneck Distance \cite{cohen2007stability}] \label{def:Hausdorff and bottleneck}
The 
bottleneck distance between two point sets $\mathcal{Y}$ and $\mathcal{Z}$ is given by
\begin{align}
%
d_B(\mathcal{Y},\mathcal{Z}) &=  \inf_{\gamma} \sup_{y} \left\lVert{y-\gamma(y)} \right\rVert_\infty ~,
\end{align}
where $y \in \mathcal{Y}$ and $z \in \mathcal{Z}$ range over all points and $\gamma$ ranges over all bijections from $\mathcal{Y}$ to $\mathcal{Z}$. Here, we interpret each point with multiplicity $k$ as $k$ individual points and the bijection is between the resulting sets.
\end{definition}

Next, we present 
triangulable spaces. As in \cite{edelsbrunner1997triangulating}, an \textit{underlying space} of a simplicial complex $K$ is defined as $\bigcup K=\bigcup_{\sigma \in K}\sigma$. If there exists a simplicial complex $K$  such that $\bigcup K$  is homeomorphic to $X$, then $X$  is triangulable. 
 
In  Fig.~\ref{fig:SPD}, we show that it two point clouds are similar, then their associated persistence diagrams are similar. That is, persistence diagrams are `stable' with respect to small changes, which can be quantified through the bottleneck distance.



\begin{theorem}[Stability Theorem \cite{bukkuri2021applications} ]\label{thm:global_stability}
 Let $X$ be a triangulable space with continuous tame functions $g, h : X \longrightarrow \mathbb{R}$, and let $D_g$ and $D_h$ denote their associated persistence diagrams obtained using a height filtration. Then the bottleneck distance between these persistence diagrams satisfies a global uniform bound
 $$d_B(D_g,D_h)\leq \left\lVert {g-h} \right\rVert_\infty .$$
\end{theorem}

While Thm.~\ref{thm:global_stability} describes the stability of persistence diagrams with respect to changes to the functions to which filtrations are applied, one can also use   it to equivalently prove the stability of persistence diagrams for when the function changes due to perturbations of points. We present the following example.

\begin{corollary}[Stability of VR Persistence Diagrams to Point Perturbations]\label{cor:point_stability}
Let $\mathcal{Y}=\{y_i\}_{i \in \mathcal{V}}\in\mathbb{R}^d$ with $\mathcal{V}=\{1,...,N\}$ denote a set of $N$ points, and let $\mathcal{Y'}=\{y_i'\}_{i \in \mathcal{V}}\in\mathbb{R}^d$ be a set of perturbed points
such that $||y_i-y_i'||_2 \le \epsilon$ for all $i$. 
If
$g: \mathcal{Y} \longrightarrow \mathbb{R}^{n \times n}_+$
denotes the pairwise distance function in which $g_{ij}(\mathcal{Y})=||y_i-y_j||_2$ and $h: \mathcal{Y'} \longrightarrow \mathbb{R}^{n \times n}_+$
 with $h_{ij}(\mathcal{Y'})= ||y'_i-y'_j||_2$. It then follows that
$d_B(D_g,D_h)  \le 2\epsilon$.

\begin{proof}
\begin{align}
d_B(D_g,D_h) \le ||g-h||_\infty&=\max_{i,j\in\mathcal{V}} |g(y_i,y_j)-h(y'_i,y'_j)|\nonumber\\
&=\max_{i,j\in\mathcal{V}} \left|   ||y_i-y_j||_2-||y'_i-y'_j||_2 \right|\nonumber .
\end{align}
However, if we define  $y'_i = y_i + e_i$, then $||e_i||_2\le\epsilon$ for any $i\in\mathcal{V}$ and
\begin{align}
||y'_i-y'_j||_2 = ||(y_i-y_j) + (e_i-e_j) ||_2 \le ||(y_i-y_j)||_2 + || (e_i-e_j) ||_2.\nonumber
\end{align}
It follows that
\begin{align}
\left|   ||y_i-y_j||_2-||y'_i-y'_j||_2 \right| \le || (e_i-e_j) ||_2 \le 2\epsilon \nonumber ,
\end{align}
which subsequently also bounds  $d_B(D_g,D_h) $.

\end{proof}
\end{corollary}

\subsubsection{Convergence}\label{sec:convergence}

The stability of persistence diagrams also has important consequences for the convergence of a sequence of persistence diagrams that is associated with   convergent sequence of point clouds.
Focusing on VR filtrations, consider a sequence  of point clouds  $\mathcal{Y}^m = \{y^m_i\}_{i=1}^N  $ of fixed size $N = |\mathcal{Y}^m |$ for each $ m \in \mathbb{N}_{+}$ with   $y_i^m \in \mathbb{R}^p$. Assume for each point $i$ that the sequence converges  $y_i^m \to y_i$  such that $\lVert{y^m_i-y_i}\rVert_2  \le 1/m$. 
Note that such a convergence criterion can be ensured by considering a subsequence  in which each subsequent element is chosen so that the bound is true for all $i\in\{1,\dots,N\}$.
Let $D_{f_m}$ and $D_f$ be the persistent diagrams for VR filtrations using the pairwise distance functions $f_m$ and $f$, which are respectively defined using the respective point clouds $\mathcal{Y}^m$ and $\mathcal{Y}$.
The Stability Theorem and Corollary~\ref{cor:point_stability} imply convergence of the associated persistence diagrams since
$$ d_B(D_{f_m},D_f)\leq \left\lVert{f_m-f} \right\rVert_\infty \le 2 \lVert{y^m_i-y_i}\rVert_2 \le   \frac{2}{m} ,$$ which converges to $0$ as $m\to\infty$.

\section{Topological data analysis using kNN orderings}\label{sec:main}

We now present our main findings: an approach for  persistent homology that is  based on the relative positioning of points according to their $k$-nearest-neighbor (kNN) sets. Our approach
introduces kNN complexes and filtrations using a filtration parameter $k\in\mathbb{N}_+$, thereby contrasting  VR and other filtrations that use a distance threshold $\epsilon$ as the filtration parameter. We will show that kNN-based persistent homology has certain advantages  that can benefit applications for which the relative positioning of data points is important, i.e., as opposed to their precise locations.

This section is organized as follows. In Sec.~\ref{sec:NNO}, we define kNN  orderings of points. In Sec.~\ref{sec:NNO_ph}, we develop kNN filtrations, use them to study kNN persistent homology, and compare them to the study of VR filtrations. In Sec.~\ref{sec:NNO_stab}, we analyze the stability and convergence of persistence diagrams  resulting from kNN filtrations.

\subsection{kNN orderings, neighborhoods,  and symmetrization}\label{sec:NNO}

We begin with a definition for kNN orderings. One complication is that the ordering of nearest neighbors is not necessarily symmetric, and so we will also describe additional transformations that we will use to symmetrize orderings prior to    constructing filtrations.

\begin{definition}[$k$-Nearest-Neighbor Orderings] \label{def:NNO}
Let $\mathcal{V} = \{1,\dots,N\}$ enumerate a set of points $\mathcal{Y}=\{y_i\}_{i\in\mathcal{V}}\in\mathbb{R}^p$ in a normed metric space and $\{f_{ij} = || y_i-y_j||_2 \}$ be their pairwise distances resulting from the pairwise distance map given in Definition~\ref{def:NNOF}. For each $i$, let  $k_{ij}$ denote the \textbf{nearest-neighbor order} of $y_j$ with respect to $y_i$ so that $y_j$ is the ($k_{ij}$)-th nearest neighbor of  $y_i$ ($k_{ii}\equiv 0$). The ordering is formally defined for each fixed $i$ by  $\{k_{ij}\}_{j=1}^N = \phi(\{f_{ij}\}_{j=1}^N)$, where
$\phi: \mathbb{R}^{N} \longrightarrow \mathbb{N}^{N}$ is the ``argsort function''.
\end{definition}

In the above, we assume that the orderings are well-defined in that for a given $i\in\mathcal{V}$, the entries in the set $\{f_{ij}\}_{j=1}^N$ are unique. In practice, if two distances $f_{ij}$ and $f_{ij'}$ are the same and they correspond or order $k$, then we assign an ordering at random so that either $k_{ij} = k$ and $k_{ij'} = k+1$, or vice versa. If even more distances are the same, then we similarly assign them a relative ordering uniformly at random. In principle, the argsort function $\phi$ can be defined in a variety ways to handle the ordering of repeated entries. Given the above definition of kNN orderings, we now introduce an associated  map between a set of points and an associated matrix with entries $k_{ij}$.

\begin{definition}[kNN Ordering Function] \label{def:NOf}
Let $\mathcal{Y}=\{y_i\}_{i \in \mathcal{V}}$ denote a set of $N\ge 2$ points with $\mathcal{V}=\{1,2,...,N\}$ in a Euclidean metric space and recall the row-defined pairwise distance function $f_i: \mathcal{Y} \longrightarrow \mathbb{R}^{N }_{+}$ given in Definition~\ref{def:NNOF}. 
Letting $\phi$ be the argsort function, we define the map
\textbf{k-nearest-neighbor (kNN) ordering function} $F:\mathcal{Y}\longrightarrow \mathbb{N}^{N\times N }_{+}$ as the matrix-valued function $F=[F_1,\dots,F_N]^T$ in which each row $F_i$ is defined by $F_i = \phi \circ f_i: \mathcal{Y} \longrightarrow \mathbb{N}^{N }_{+}$.
\end{definition}

We will use kNN orderings of points to define   local neighborhoods $\mathcal{N}_{ik}$ for each point $i\in\mathcal{V}$, noting that any nested sequence of neighborhoods defines a local filtration of the vertices $\mathcal{V}$. 

\begin{definition}[kNN Neighborhoods] \label{def:local_NNO_filtration}
Given a set of kNN orderings $\{k_{ij}\}$ for $i,j \in\mathcal{V}$, we define the \textbf{kNN neighborhoods} as the sets $\mathcal{N}_{ik} = \{j  | k_{ij} \le k \}$ for each  $i\in\mathcal{V}$ and $k \in \{0\}\cup \mathcal{V}\setminus \{N\}$. 
\end{definition}

Importantly,  the kNN orderings and kNN neighborhoods are not symmetric---that is,   $j\in\mathcal{N}_{ik}$   does not imply $i\in\mathcal{N}_{jk}$.  Because we would like to use the kNN neighborhood sets to construct filtrations of simplicial complexes that contain undirected simplices, we will now describe how to construct symmetric kNN neighborhoods.

\begin{definition}[Symmetrized kNN Orderings and Neighborhoods] \label{def:local_NNO_symm}
We define three types of \textbf{symmetrized kNN  neighborhoods}---${\mathcal{N}}_{ik}^{min}$, ${\mathcal{N}}_{ik}^{trans}$ and ${\mathcal{N}}_{ik}^{max}$---that use, respectively, the following three symmetrizations of kNN orderings  $\{k_{ij}\}$:
\begin{align}\label{eq:sym}
    \tilde{k}_{ij} &= \min\{ k_{ij} , k_{ji} \} \nonumber\\
    \overline{k}_{ij} &= (k_{ij} + k_{ji} )/2 \nonumber\\
     \hat{k}_{ij} &= \max\{ k_{ij} , k_{ji} \}.
\end{align}
\end{definition}
Notably, these kNN neighborhoods  satisfy the following nestedness relation.
\begin{lemma}[Nestedness of Symmetrized kNN Neighborhoods] \label{lemma:monotonic2}
Consider fixed $k$ and $i$. Then the neighborhood sets $ \mathcal{N}_{ik}^{min}, \mathcal{N}_{ik}^{trans},  \mathcal{N}_{ik}^{max} $  satisfy the following nestedness relationships:
\[
\mathcal{N}^{max}_{ik} \subseteq \mathcal{N}^{trans}_{ik} \subseteq \mathcal{N}^{min}_{ik}.
\]

\begin{proof}
For any $k_{ij},k_{ji}\in\mathbb{R}$, one has $\tilde{k}_{ij} \le \overline{k}_{ij} \le \hat{k}_{ij}$ by definition of the min and max functions.
If $j \in \mathcal{N}^{max}_{ik}$, then $\hat{k}_{ij} \leq k$. This implies $\tilde{k}_{ij} \le \overline{k}_{ij} \le \hat{k}_{ij} \le k$, and so $j \in \mathcal{N}^{min}_{ik} \subseteq \mathcal{N}^{trans}_{ik} \subseteq \mathcal{N}^{max}_{ik}$.
\end{proof}
\end{lemma}

We now   define graphs and clique complexes using symmetrized kNN orderings and neighborhoods.

\begin{definition}[Symmetrized kNN Graph] 
\label{def:NNO_graph}
Given a set of points enumerated by  $\mathcal{V}$, let $\mathcal{E}_k = \{(i,j)| i\in \mathcal{V}, j\in \mathcal{N}_{ik}\}$ be a set of undirected edges connecting pairs of nearest neighbors, which are symmetrically defined by choosing $\mathcal{N}_{ik}\in\{ \mathcal{N}_{ik}^{min}, \mathcal{N}_{ik}^{trans},  \mathcal{N}_{ik}^{max} \}$. We denote the respective graphs $\mathcal{G}_k(\mathcal{V},\mathcal{E}_k)$, depending on symmetrization method, by
$\mathcal{G}_{k}^{min}$, $\mathcal{G}_{k}^{trans}$,  and $\mathcal{G}_{k}^{max}$.
\end{definition}

\begin{definition}[Symmetrized kNN Clique Complex] \label{def:NNO_global}
Given a symmetrized kNN graph $\mathcal{G}_k \in\{\mathcal{G}_{k}^{min}, \mathcal{G}_{k}^{trans},  \mathcal{G}_{k}^{max}\}$ with vertices $\mathcal{V}$ and edges $\mathcal{E}_k$, we construct its associated clique complex $X_k = Cl(\mathcal{G}_{k} )$, and we similarly use a superscript to denote the method of symmetrization, i.e., $X_k^{min}$, $X_k^{trans}$, and $X_k^{max}$. 
\end{definition}

\begin{corollary}
[Nestedness of Symmetrized kNN Graphs and kNN Complexes] \label{lemma:monotonic3}
Symmetrized kNN graphs and kNN complexes satisfy the following nestedness relations
\begin{align}
\mathcal{G}^{max}_{k} &\hookrightarrow \mathcal{G}^{trans}_{k} \hookrightarrow \mathcal{G}^{min}_{k}\nonumber\\
X^{max}_{k} &\hookrightarrow X^{trans}_{k} \hookrightarrow X^{min}_{k}.
\end{align}
\begin{proof}
The results follow immediately from  Lemma~\ref{lemma:monotonic2}, which proved the nestedness of kNN neighborhoods.
\end{proof}
\end{corollary}

Corollary~\ref{lemma:monotonic3} implies that $\mathcal{G}^{max}_{k}$ is a subgraph of $\mathcal{G}^{trans}_{k}$, which is a subgraph of $\mathcal{G}^{min}_{k}$. Similarly, $X^{max}_{k}$ is a subcomplex of $X^{trans}_{k}$, which is a subcomplex of $X^{min}_{k}$.

\subsection{Filtrations and persistent   homology using kNN complexes}\label{sec:NNO_ph}

We now formulate filtrations and persistent homology using symmetrized kNN complexes. Varying $k$ gives rise to a sequence of nested sets called a filtration, and we will define several types based on the different symmetrization methods.
We will also define local and global versions of filtrations.

\begin{figure}[htp]
\begin{center}
    \includegraphics[width=\linewidth]{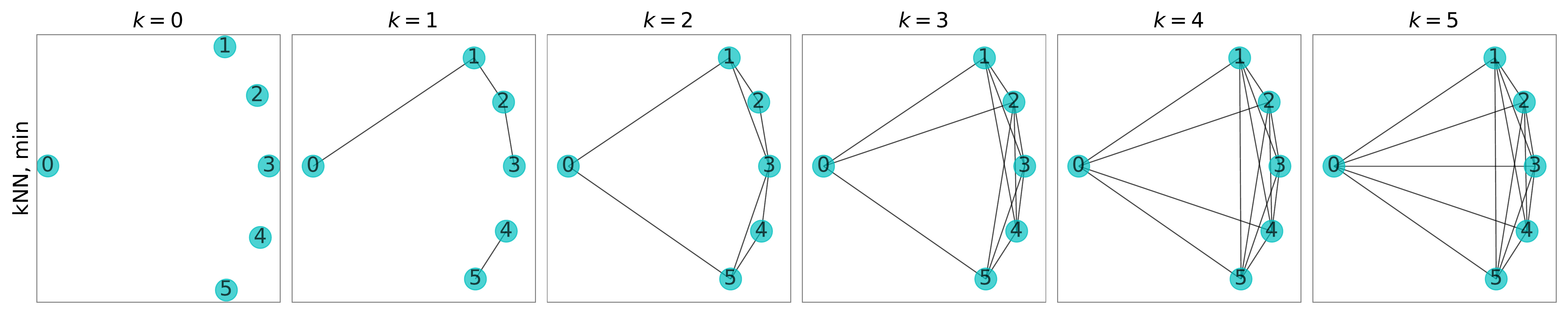}
    \includegraphics[width=\linewidth]{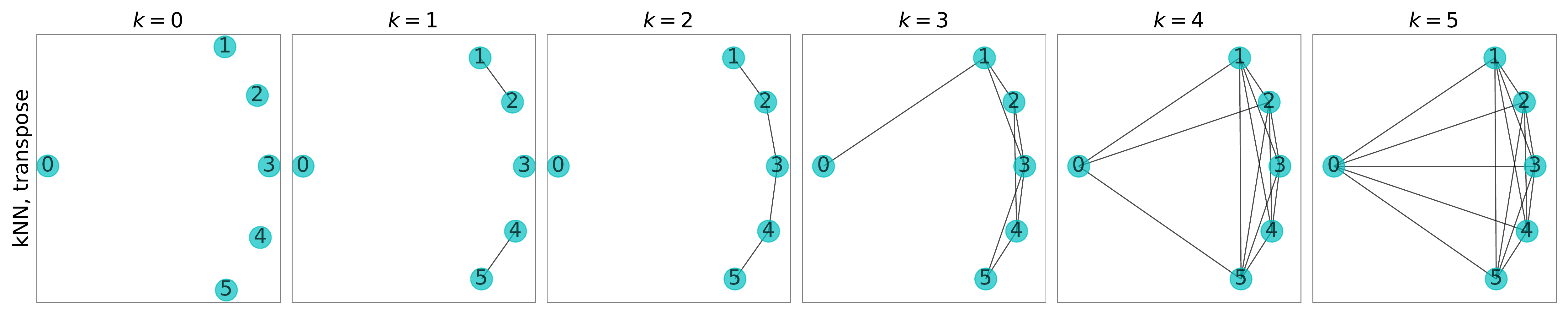}
    \includegraphics[width=\linewidth]{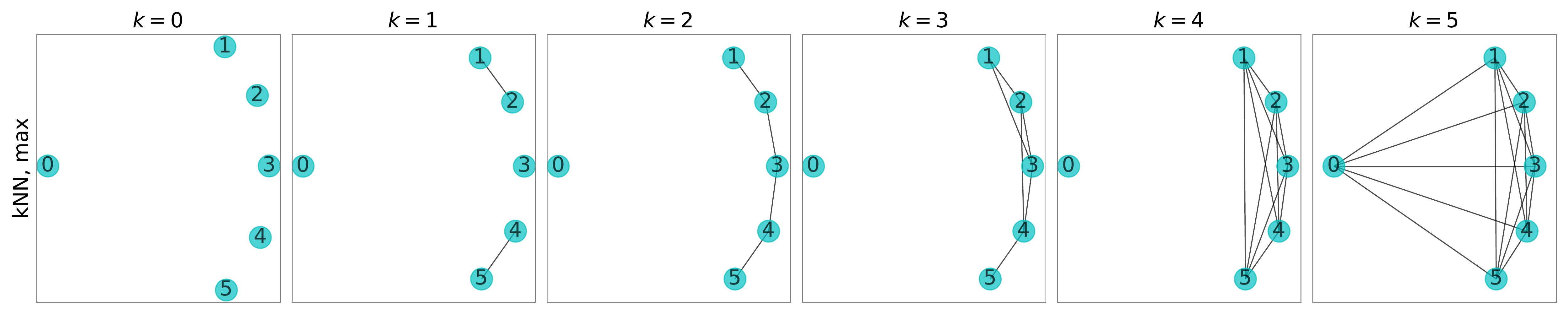}
    \caption{Visualization of  kNN-filtered simplicial complexes for a point cloud using the three types of symmetrization given by Definition~\ref{def:local_NNO_symm}.
    Comparing across the columns, observe the nestedness property given by Corollary~\ref{lemma:monotonic3}.
    }\label{Fig:NNOs_sym}
\end{center}
\end{figure}

\begin{definition}[Global kNN Filtration]
\label{def:NNO_globa}
Let $X_k,X_{k'}\in \{ X_k^{min},X_k^{trans},X_k^{max}\}$ be   kNN complexes with the same symmetrization method.
Then we define a \textbf{kNN-filtered simplicial complex}
\[
\{X_k \longrightarrow X_{k'} \}_{0\leq k \leq k^\prime\le N-1} .
\]
\end{definition}

In Fig.~\ref{Fig:NNOs_sym},  we illustrate global kNN filtrations with the   three different symmetrization methods given by Definition~\ref{def:local_NNO_symm}. For simplicity, we only visualize 0-simplices and 1-simplices. By comparing the kNN complexes across a given column, one can observe the nestedness relations defined by Corollary~\ref{lemma:monotonic3}.

\begin{figure}[t]
\begin{center}
     \includegraphics[width=.45\linewidth]{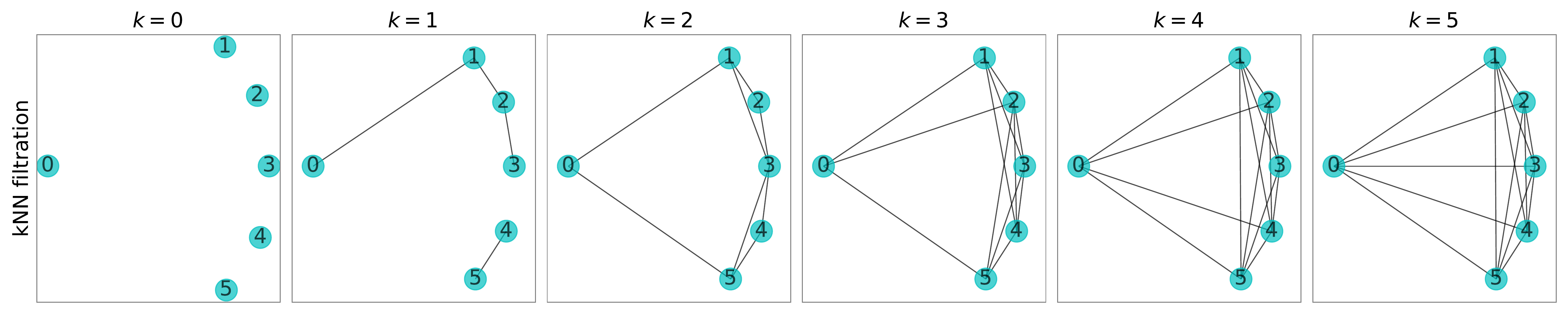}~~~
     \includegraphics[width=.45\linewidth]{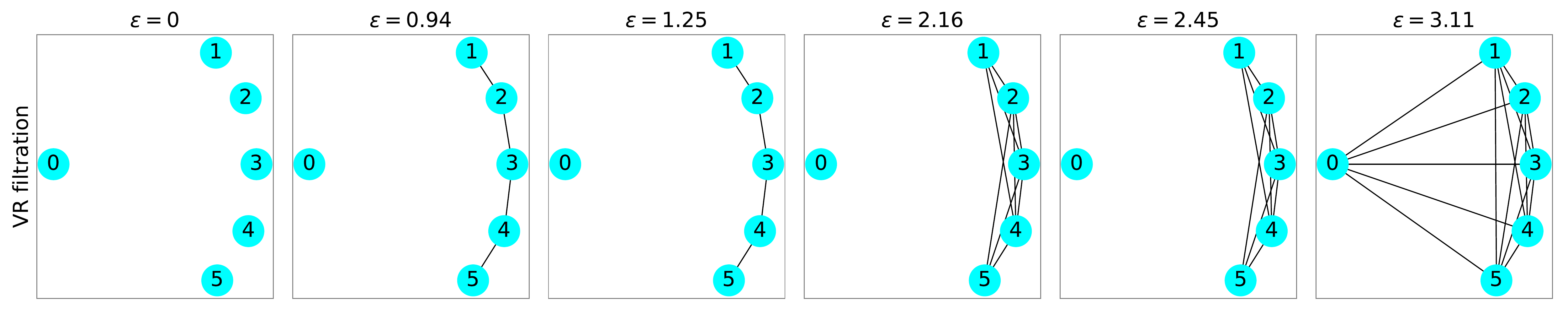}
     \includegraphics[width=.45\linewidth]{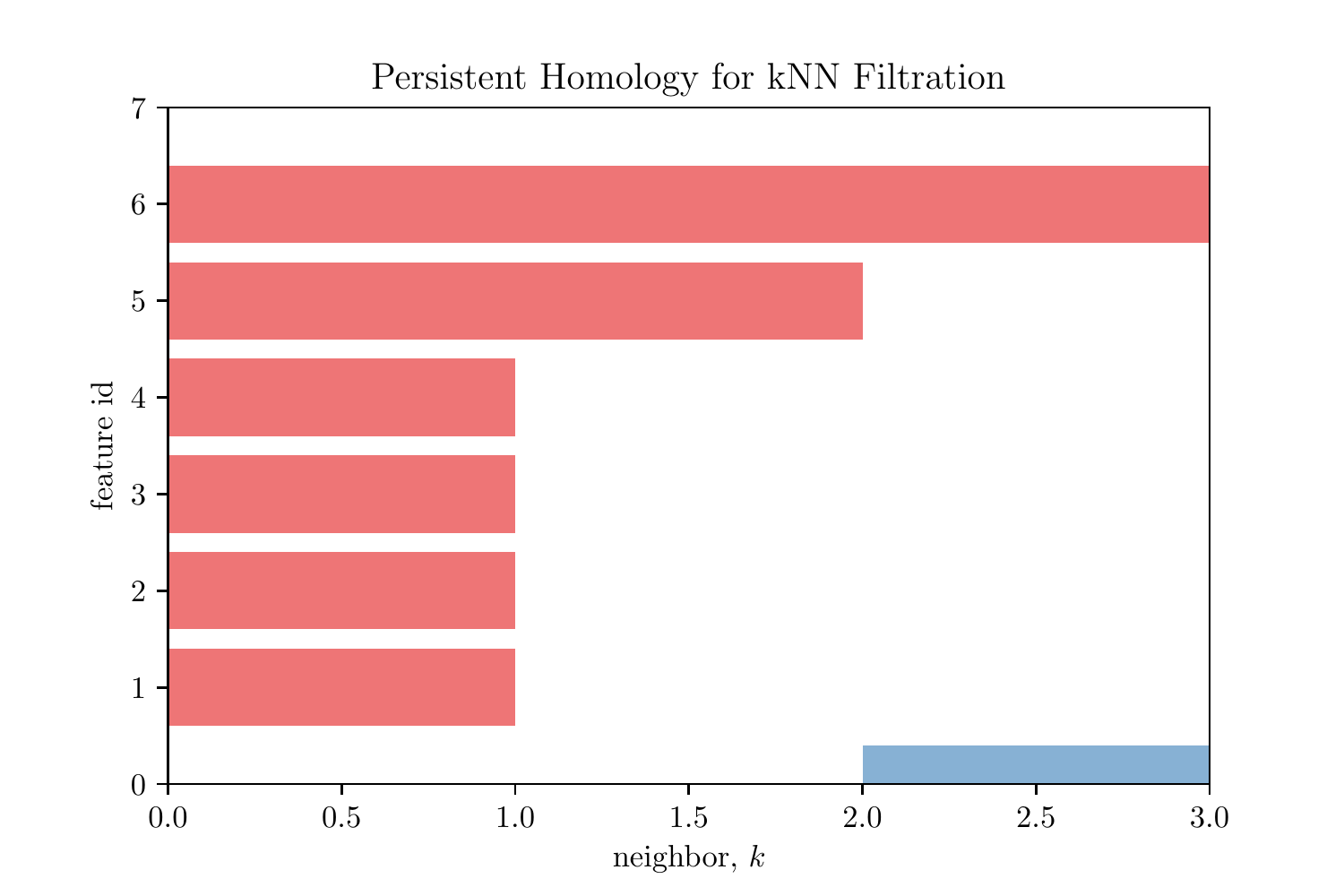}
    \includegraphics[width=.45\linewidth]{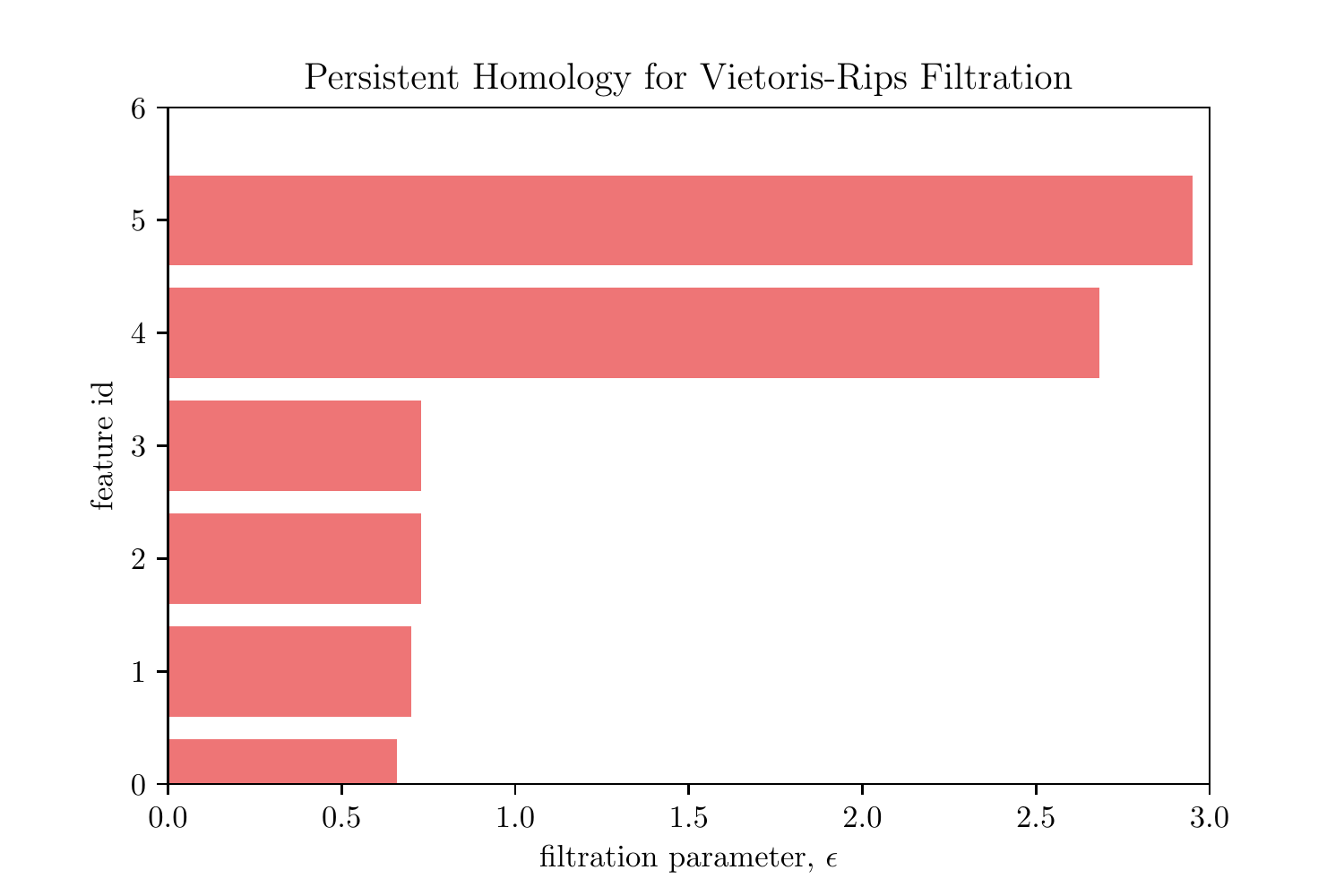}
     \caption{Comparison of persistent homology for an example point cloud using two different filtrations: (left) our proposed kNN filtration; and (right) a VR filtration. The red and blue persistence barcodes indicate 0-dimensional and 1-dimensional cycles, respectively.
     Observe that the kNN filtration reveals a 1-cycle, whereas the VR filtration does not.
     }\label{Fig:VR_filtration}
     \end{center}

\end{figure}

We formulate  persistent homology for kNN filtrations of a simplicial complex anaologous to that defined for a VR filtrations (recall Sec.~\ref{sec:holomology}), and in Fig.~\ref{Fig:VR_filtration}, we visualize  persistence diagrams for both (left) a  kNN filtration and (right) a VR filtration for an example point cloud. The red and blue persistence barcodes indicate 0-dimensional and 1-dimensional cycles, respectively. Observe that the kNN filtration reveals a 1-cycle, whereas the VR filtration does not. Hence, filtrations of simplicial complexes (specifically clique complexes) according to pairwise distances and kNN provide complementary homological information.

In the next sections, we will further compare persistence diagrams resulting from kNN and VR filtrations;
however, it's worth highlighting several differences here. 
First, one  benefit of using kNN sets versus a distance threshold $\epsilon$ is that filtrations are standardized. That is, the filtration parameter range $k\in[0,N]$ for kNN filtrations is always the same for a set of $N$ points. In contrast, different point clouds have different lengths scales, and so the range of a distance-based filtration parameter is in general not standardized.  One could seek to standardize the ranges of VR filtration parameters by standardizing distances; however, there are many ways to normalize a point cloud, such as dividing by the mean distance or by the distances' standard deviation. Given that one could implement a variety of normalization approaches, comparing VR filtrations across different point clouds that have different length scales or different dimensions is not straightforward. In contrast, kNN filtrations are standardized by definition.

Second, the filtration parameter $\epsilon\in\mathbb{R}_+$ for VR filtrations can in principle take on any positive number. In contrast, $k\in\{0,\dots, N-1\}$ can only take on integer values (or in the case of the symmetrization method of $trans$, half integers as well). Potentially, this reduced space of possible filtration parameters  could benefit the computationally efficiency of implementing kNN filtrations, although we don't explore that pursuit herein. That said, our experiments generally find kNN filtrations to have higher computational complexity, since they require both computing pairwise distances as well as their orderings.

\subsection{Stability and convergence of kNN homology}\label{sec:NNO_stab}

Here,  we will study two key properties for the space of persistence diagrams that follow from kNN filtrations: stability and convergence.  


\subsubsection{kNN-preserving transformations}\label{sec:NNO-preserve}

We observe that if points undergo  perturbations that are sufficiently small, then it is possible to move points without changing any of the neighbor orderings. To make this more precise, we define a family of point-cloud transformations that have this property.

\begin{definition}[Local and Global  kNN-Preserving Transformations] \label{def:NNO_transformation}
Let $\mathcal{Y}=\{y_i\}_{i \in \mathcal{V}}$ with $ y_i \in \mathbb{R}^p$ be a set of points that are enumerated $\mathcal{V}=\{1,...,N\}$,   $f: \mathcal{Y} \longrightarrow \mathbb{R}^{N \times N}_{+}$ be the pairwise distance function given by Def.~\ref{def:NNOF}, and $F: \mathcal{Y} \longrightarrow \mathbb{N}^{N \times N}_{+}$ be the neighbor-ordering function given by Def.~\ref{def:NOf}. 
We define a transformation $h: \mathcal{Y} \longrightarrow \mathcal{Y'}$ with $h(x)=(h_1(y_1),...,h_n(y_n))=(y'_1,...,y'_n)$ and define two types of preservation properties:

\begin{itemize}
    \item $h$  is local kNN-preserving  for some $i \in \mathcal{V}$  if and only if $(F \circ h)_{ij}(\mathcal{Y}) = F_{ij}(\mathcal{Y})$  for fixed $i$ and $j \in \mathcal{V}$.
    \item $h$  is global  kNN-preserving if and only if it is local kNN-preserving for all $i\in\mathcal{V}$. 
\end{itemize}
\end{definition}

Note that local kNN-preservation implies that   the neighbor orderings are preserved for a particular point $x_i$, while global-kNN preservation implies that  all neighbor orderings are preserved for all points. We also note that local and global  kNN-preserving transformations can also be similarly defined for the three symmetrized versions of kNN orderings.
We also define a slight variation in which the preservation of orderings is  only required for the   nearest neighbor orderings up to a finite size $k\le K$, which we refer to as ``$K$-bounded'' local and global kNN-preserving transformations.

As one example, consider a set of three points $x_i\in \{0,2,3\}\subset \mathbb{R}$ and the transformation $h(\{0,2,3\}) = \{0,2.1,3\} $. It is both local and global kNN-preserving, since the perturbation $ 2\to2.1$ is sufficiently small such that none of the neighbor orderings change. In fact, all isometric transformations including (i.e., translation, rotation, reflection and glide reflection) are global kNN-preserving transformations because they leave the pairwise distances unchanged.

\begin{figure}[htp]
\begin{center}
    
\includegraphics[width=1\linewidth]{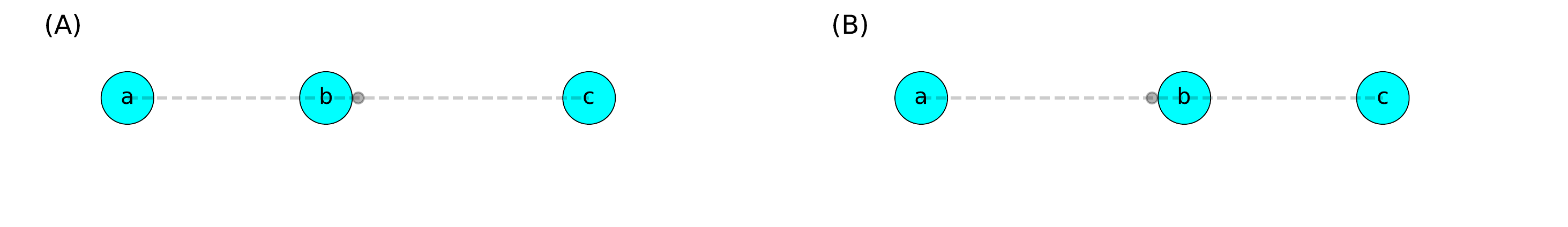}
\caption{Example with 3 points $\mathcal{Y} = \{x_a,x_b,x_c\}$ with $x_a=-1$, $x_c =1$ and either (A) $x_b=-\epsilon$ or (B) $b=\epsilon$. Note that the perturbation can be made arbitrarily small for any  $\epsilon>0$, and the nearest-neighbor orderings are different. (See Tables \ref{table:Rank line origin} and \ref{table:Rank line transform}.) }
\label{fig:stab}
\end{center}

\end{figure}

However, there are many transformations that are not  kNN preserving. In Fig.~\ref{fig:stab}, we illustrate an example where a small perturbation to one point can change the neighbor ordering in a discontinuous way. Consider a set of points $\mathcal{V}=\{a,b,c\}$ with locations $\mathcal{Y}=\{x_a,x_b,x_c\}\subset \mathbb{R}$ where $x_a=-1$, $x_c=1$ and $x_b=-\epsilon$ (for some small value $0<\epsilon \ll 1$). Then the neighbor ordering matrix is given in Table~\ref{table:Rank line origin}.
Now consider the transformation $h_1(x_a)=x_a$, $h_2(x_b)=x_b+2\epsilon$, and $h_3(x_c)=x_c$, then the neighbor ordering matrix changes, as shown in Table~\ref{table:Rank line transform}.

\begin{table}[ht]
\begin{minipage}[b]{0.4\linewidth}
\centering
{\small{\begin{tabular}{ | l | r | r | r | r |}
    \hline
     &  point a &  point b & point c\\ \hline \hline
     point a & 0 & 1&2 \\ \hline
     point b & 1 & 0&2 \\ \hline
     point c & 2 & 1&0 \\ \hline
    \end{tabular}}}
    \caption{kNN orderings before transformation $h$.}
    \label{table:Rank line origin}
\end{minipage}\hfill
\begin{minipage}[b]{0.4\linewidth}
\centering
{\small{\begin{tabular}{ | l | r | r | r | r |}
    \hline
     &  point   a &  point b &  point c \\ \hline \hline
     point a & 0 & 1 &2 \\ \hline
     point b & 2 & 0& 1 \\ \hline
     point c  & 2 & 1&0 \\ \hline
    \end{tabular}}}
    \caption{kNN orderings after transformation $h$.}
    \label{table:Rank line transform}
\end{minipage}
\end{table}

Given the definitions of local and global kNN-preserving transformations, we 
can also define equivalence classes for any finite set of points $ \mathcal{Y}= \{y_i\} \subset \mathbb{R}^p$.

\begin{definition}[kNN Equivalence of Point Sets] \label{def:NNO_equiv}
Let $F$ denote the neighbor-ordering function given by Def.~\ref{def:NOf}.
Any two point clouds $\mathcal{Y}=\{y_i\}_{i \in \mathcal{V}}$ and $\mathcal{Y'}=\{y_i'\}_{i \in \mathcal{V}}$ with $ y_i,y'_i \in \mathbb{R}^p$ are said to be kNN-equivalent if $F(\mathcal{Y}) = F(\mathcal{Y'})$. We denote the equivalence by $\mathcal{Y}\sim \mathcal{Y'}$.
Furthermore, we have  $[\mathcal{Y}] = \{ \mathcal{Y'} | \mathcal{Y'}\sim \mathcal{Y} \}$ as the kNN-equivalence class of  $\mathcal{Y}$.
\end{definition}

Indeed, one can show 

\begin{itemize}
    \item {\bf Reflexivity: $\mathcal{Y} \sim \mathcal{Y}'$.} \\
    It is trivial. Take $h(t)= \mathbf{id}$ as the kNN-preserving transformation.
    \item {\bf  Symmetry: if $\mathcal{Y} \sim \mathcal{Y'}$  then $\mathcal{Y'} \sim \mathcal{Y}$ .} \\
    Let $h: \mathcal{Y} \longrightarrow \mathcal{Y'}$ be a global kNN-preserving transformation. Then $(F \circ h)_{ij}(\mathcal{Y})=(F_{ij})(\mathcal{Y})$ and $h$ is a one-to-one function. Let $h^{-1}$ be the inverse function of $h$. So, we have: $(F \circ h^{-1})_{ij}(\mathcal{Y'})=(F_{ij})(\mathcal{Y'})$, thus $h^{-1}$ is the global kNN-preserving transformation from $\mathcal{Y'}$ to $\mathcal{Y}$.
    \item {\bf Transitivity: if $\mathcal{Y} \sim \mathcal{Y'}$ and $\mathcal{Y'} \sim \mathcal{Y''}$ then $\mathcal{Y} \sim \mathcal{Y''}$.}\\
    Let $h(t)$  be kNN-preserving transformation between $\mathcal{Y}$ and $\mathcal{Y'}$, $g(t)$ be kNN-preserving transformation between $\mathcal{Y'}$ and $\mathcal{Y''}$, then $g \circ h$ be kNN-preserving transformation between $\mathcal{Y}$  and $\mathcal{Y''}$.\\
\end{itemize}

The equivalence of kNN orderings implies that kNN filtrations and kNN-based persistent homology are identical for any two point clouds within a given equivalence class. Thus, persistence diagrams for kNN persistent homology are robust to (and in fact, unaffected by)  any kNN-preserving transformation. That said, such persistence diagrams can also be unstable to other types of perturbations, as we explore below.

\subsubsection{Stability of kNN persistence diagrams} \label{sec: stability1}

In analogy to the stability theorem for the space $\mathcal{H}_{VR}$ of persistence diagrams obtained using VR filtrations, we  consider whether the space $\mathcal{H}_{kNN}$ of persistence diagrams using kNN filtrations is also stable with respect to perturbations of the associated point cloud.
We fist show that persistence diagrams resulting from kNN filtrations are do not satisfy the same stability condition as that for VR filtrations.

\begin{theorem}[kNN   Persistence Diagrams are not Uniformly Globally Stable to Point Changes]\label{thm:kNN_sta}
Let $\mathcal{Y}=\{y_i\}_{i \in \mathcal{V}}$  with $\mathcal{V}=\{1,...,N\}$  be a point cloud with a neighbor-ordering function $F=\phi \circ f$, where
$f: \mathcal{Y} \longrightarrow \mathbb{R}^{N \times N}_{+}$ is the pairwise distance function and
$\phi: \mathbb{R}^{N \times N}_{+} \longrightarrow \mathbb{N}^{N \times N}_{+}$ is the `argsort' function. Further let   $\mathcal{Y'}=\{y'_i\}_{i \in \mathcal{V}}$ be a second point cloud with the neighbor-ordering function $G=\phi \circ g$.
Further let $D_F$ and $D_G$ denote the persistence diagrams of kNN filtrations applied to $\mathcal{Y}$ and $\mathcal{Y'}$. 
Further, define the metric for two enumerated point clouds $\left\lVert{\mathcal{Y}-\mathcal{Y'}} \right\rVert_\infty= \max_{i\in\mathcal{V}} ||y_i-y_i'||_2$.
Then there does not exist a uniform global bound on the  bottleneck distance between persistence diagrams of the form:
$$
d_B(D_F,D_G) \leq \lVert{F-G}\rVert_\infty 
=\left\lVert{\phi \circ f- \phi \circ g} \right\rVert_\infty \not\leq L \left\lVert{\mathcal{Y}-\mathcal{Y'}} \right\rVert_\infty.
$$

\begin{proof}
The first inequality is satisfied by applying the original statement of the stability theorem: persistence diagrams are uniformly bound by the maximal difference between the functions  that are filtered. By definition, the neighbor-ordering function yields a matrix of k nearest neighbors, $[k_{ij}] = F(\mathcal{Y})$ and $[k'_{ij}] = G(\mathcal{Y'})$, implying 
$$
\lVert{F-G}\rVert_\infty \equiv \max_{ij} |k_{ij}-k'_{ij}|.
$$
That is, the difference in persistence diagrams is bound by the maximum difference in neighbor orderings.  However,  the stability of persistent homology with respect to neighbor-ordering changes does not imply stability with respect to point-location changes. 

To disprove the last inequality, we provide a counter example. Recall the sets of $N=3$ points in Fig.~\ref{fig:stab} with  $\mathcal{V}=\mathcal{V}'=\{a,b,c\}$ and locations $\mathcal{Y} = \{1,-\epsilon,1\} \in\mathbb{R}$ and $\mathcal{Y}' = \{1,\epsilon,1\} \in\mathbb{R}$.  By construction, $||F-G||_\infty=1$ for any small $\epsilon$, but  $\left\lVert{\mathcal{Y}-\mathcal{Y'}} \right\rVert_\infty=2\epsilon$. Suppose there did exist a uniform bound with Lipschitz constant  $L>1$, then the points $\mathcal{Y}$ and $\mathcal{Y'}$ with any $\epsilon<1/2L$ yields a contradiction since $||F-G||_\infty=1$ and $||\mathcal{Y}-\mathcal{Y'}||_\infty<1$.

\end{proof}
\end{theorem}

\subsubsection{Topological convergence of kNN orderings} \label{sec:convergence1}

In this section, we will use the neighbor-ordering function to define a discrete-topological notion of convergence for a sequence of point sets using kNN persistent homology.

\begin{definition}
Consider a sequence of point sets $\mathcal{Y}^{(t)}=\{y^{(t)}_i| i \in \mathcal{V}\}$ with $t\in \mathbb{N}$ and $ \mathcal{V}= \{1,...,N\}$ and $y^{(t)}_i \in \mathbb{R}^p$. We say that the sequence $\{\mathcal{Y}^{(t)}\}$ has {\bf global convergence in kNN topology}  to the limit $\mathcal{Y}$  iff $\lim_{t\to \infty} F(\mathcal{Y}^{(t)}) = F(\mathcal{Y})$, where $F$ is the neighbor-ordering function. More precisely, for any $\epsilon$, there exists a $t^*$ such that $\max_{ij}|F_{ij}(\mathcal{Y}) - F_{ij}(\mathcal{Y}^{(t)})|<\epsilon$ for all $t>t^*$.
\end{definition}

\begin{definition}
Consider a sequence of point sets $\mathcal{Y}^{(t)}=\{y^{(t)}_i| i \in \mathcal{V}\}$ with $t\in \mathbb{N}$ and $\mathcal{V}= \{1,...,N\}$ and $y^{(t)}_i \in \mathbb{R}^p$. We say that the sequence $\{\mathcal{Y}^{(t)}\}$ has {\bf $\kappa$-bounded convergence in kNN topology} to  a limit $\mathcal{Y}$  iff $\lim_{t\to \infty} F_{ij}(\mathcal{Y}^{(t)})=F_{ij}(\mathcal{Y})$
for all $i \in \mathcal{V}$ and $j \in \mathcal{N}_{ik}$, where $k\le \kappa$ and $\mathcal{N}_k \in\{\mathcal{N}_{ik}^{min}, \mathcal{N}_{ik}^{trans},  \mathcal{N}_{ik}^{max}\}$ are  symmetrized neighborhoods that are defined  in Def.~\ref{def:local_NNO_filtration}.
More precisely, for any $\epsilon$, there exists a $t^*$ such that $\max_{i\in\mathcal{V},j\in \mathcal{N}_{ik}}|F_{ij}(\mathcal{Y}) - F_{ij}(\mathcal{Y}^{(t)})|<\epsilon$ for all $t>t^*$.
\end{definition}

\begin{remark}
Note that $\kappa$-bounded convergence is inclusive so that convergence for a given $\kappa$ also implies convergence for any $\kappa'\le \kappa$. In particular, global convergence in kNN topology implies  $\kappa$-bounded convergence for any $\kappa$.
\end{remark}

We now define local variants of kNN topological convergence.
\begin{definition}
Consider a sequence of point sets $\mathcal{Y}^{(t)}=\{y^{(t)}_i| i \in \mathcal{V}\}$ with $t\in \mathbb{N}$ and $ \mathcal{V}= \{1,...,N\}$ and $y^{(t)}_i \in \mathbb{R}^p$. We say that the sequence $\{\mathcal{Y}^{(t)}\}$ has {\bf $\mathcal{U}$-local convergence in kNN topology} to  a limit $\mathcal{Y}$  iff $\lim_{t\to \infty} F_{ij}(\mathcal{Y}^{(t)})=F_{ij}(\mathcal{Y})$
for all $i,j \in \mathcal{U}\subseteq \mathcal{V}$ .
More precisely, for any $\epsilon$, there exists a $t^*$ such that $\max_{i,j\in\mathcal{U}}|F_{ij}(\mathcal{Y}) - F_{ij}(\mathcal{Y}^{(t)})|<\epsilon$ for all $t>t^*$.
\end{definition}

\begin{definition}
Consider a sequence of point sets $\mathcal{Y}^{(t)}=\{y^{(t)}_i| i \in \mathcal{V}\}$ with $t\in \mathbb{N}$ and $  \mathcal{V}= \{1,...,N\}$ and $y^{(t)}_i \in \mathbb{R}^p$. We say that the sequence $\{\mathcal{Y}^{(t)}\}$ has {\bf $\kappa$-bounded $\mathcal{U}$-local convergence in kNN topology} to  a limit $\mathcal{Y}$  iff $\lim_{t\to \infty} F_{ij}(\mathcal{Y}^{(t)})=F_{ij}(\mathcal{Y})$
for all $i \in \mathcal{U}\subseteq \mathcal{V}$ and $j \in \mathcal{N}_{ik}$, where the  symmetrized neighborhoods 
$\mathcal{N}_k \in\{\mathcal{N}_{ik}^{min}, \mathcal{N}_{ik}^{trans},  \mathcal{N}_{ik}^{max}\}$
are given  in Def.~\ref{def:local_NNO_filtration} and $k\le \kappa$.
More precisely, for any $\epsilon$, there exists a $t^*$ such that $\max_{i\in\mathcal{U},j\in \mathcal{V}}|F_{ij}(\mathcal{Y}) - F_{ij}(\mathcal{Y}^{(t)})|<\epsilon$ for all $t>t^*$.
\end{definition}

Note   for any point set that convergences in global kNN topology (or k-bounded convergence in kNN topology),  that there exists a $t^*$ such that $F_{ij}(\mathcal{Y}) = F_{ij}(\mathcal{Y}^{(t)})$ for all $t>t^*$ and $i,j \in\mathcal{V}$ (or $i \in\mathcal{V}$ and $j\in N_{ik}$).
Specifically, we can choose $\epsilon\in(0,1)$. Since $F_{ij}(\mathcal{Y})\in\mathbb{N}$ for any point set $\mathcal{Y}$, the condition $|F_{ij}(\mathcal{Y})-F_{ij}(\mathcal{Y}^{(t)})|<\epsilon$ implies $F_{ij}(\mathcal{Y})-F_{ij}(\mathcal{Y}^{(t)})=0$ for all $t>t^*$.
That is, kNN convergence is a discrete property that is exactly obtained for sufficiently large $t>t^*$. 
This significantly contrasts the notion of convergence in norm, which is often asymptotically approached rather than exactly obtained.
The next theorem more precisely establishes a relation between convergence in a normed metric space and convergence in kNN topology.

\begin{theorem}[Convergence in Norm  Implies Convergence in kNN Topology]\label{thm: convergence_NNO}
Let $\mathcal{Y}^{(t)} = \{y^{(t)}_i| i \in \mathcal{V}\}\subset\mathbb{R}^p$ be an element of a sequence of point clouds in which each point converges in norm to some limit,
$\lim_{t\to\infty}||y_i^{(t)} - y_i||=0$ for each $i\in \mathcal{V}$. 
%
Define $\mathcal{Y} = \{y_i | i\in\mathcal{V} \} \subset \mathbb{R}^p$ as that limit.
Then $\mathcal{Y}^{(t)} \to \mathcal{Y}$ in kNN topology, and there exists a $t^*$ such that $F(\mathcal{Y}^{(t)})=F(\mathcal{Y})$ for any $t>t^*$, where $F$ is the neighbor-ordering given by Def.~\ref{def:NOf}.

\begin{proof}
To prove this theorem, we use a contradiction. 
If $\mathcal{Y}^{(t)}$ does not converge to $\mathcal{Y}$ in kNN-topology, then for any $t$, there exists at least one neighbor ordering that differs, i.e., $F_{ij}(\mathcal{Y}^{(t)}) \neq F_{ij}(\mathcal{Y})$ for some $i,j \in\mathcal{V}$. Letting  $d_{ij}^{(t)} = ||y_i^{(t)} - y_j^{(t)}||_2 $ and  $d_{ij}^{(\infty)} = ||y_i  - y_j ||_2 $, this implies   $d_{ij}^{(\infty)} - d_{ik}^{(\infty)}>0$ but $d_{ij}^{(t)} - d_{ik}^{(t)} \le 0$ for some $i,j, k \in\mathcal{V}$. (This must be true of the distances if the neighbor orderings differ.) Let $\delta^*= \max_{i,j,k} |d_{ij}^{(\infty)} - d_{ik}^{(\infty)}|$. However, because $\lim_{ t\to\infty} || y_i^{(t)} -y_i|| $ for each $i$, there exists a $t^*$ such that $|| y_i^{(t)} -y_i||_\infty < \epsilon/4$ for all $i$ when $t >t^*$. It then follows that $ |d_{ij}^{(t)} - d_{ij}^{(\infty)}| < \epsilon/2$ and $|d_{ik}^{(t)} - d_{ik}^{(\infty)}| < \epsilon/2$.
Now suppose  $d_{ij}^{(\infty)} - d_{ik}^{(\infty)} >0$.
Then   
\begin{align}
d_{ij}^{(t)}-d_{ik}^{(t)} &> d_{ij}^{(t)} - (d_{ik}^{(\infty)} + \epsilon/2) \nonumber\\
&> (d_{ij}^{(\infty)} - \epsilon/2) -  (d_{ik}^{(\infty)} + \epsilon/2) \nonumber\\
&=  d_{ij}^{(\infty)} - d_{ik}^{(\infty)} - \epsilon\nonumber\\
&=  \delta^* - \epsilon.
\end{align}
This is true for any $\epsilon$, so we may choose $\epsilon<\delta^*$, which implies $d_{ij}^{(t)}-d_{ik}^{(t)}>0$, contradicting the above statement. 

\end{proof}
\end{theorem}

Finally, we prove that the reverse property does not necessarily hold.

\begin{theorem}[Convergence in kNN Topology Does Not Imply Convergence in Norm] \label{thm: convergence NNO1}

Let $\mathcal{Y}^{(t)} = \{y^{(t)}_i| i \in \mathcal{V} \}\subset\mathbb{R}^p$ be an element of a sequence of point sets that converges in kNN topology.
Then each $y_i^{(t)}$ may or may not converge to some limit $y_i$ as $t\to \infty$.
 \begin{proof}
See the proof to Thm.~\ref{thm:kNN_sta} for a counter example in which points converge as $\epsilon\to 0$, but not for kNN topology, which remains bounded below by 1.

\end{proof}
\end{theorem}

\section{kNN persistent homology reveals topological convergence of pageRank algorithm} \label{sec:PagerankandRD}

In this section, we   study the convergence of an iterative method for approximating Google's PageRank an use   kNN   persistent homology to develop a perspective from discrete topology---that is, as opposed to the typical geometric perspective of convergence under a vector norm.
%
In Sec.~\ref{sec:PR1}, we present the PageRank algorithm.
In Sec.~\ref{sec:Dolph}, we    present numerical experiments for an empirical network: a dolphin social network.
In Sec.~\ref{sec:Dolpu}, we study PageRank for the same network using   local versions of kNN topological convergence.

\subsection{PageRank}
\label{sec:PR1}
Google developed \emph{PageRank} to solve the problem of web search and ranking for the World Wide Web \cite{brin1998anatomy}. Their aim was to create an importance measure for each webpage distinguish highly recognizable, relevant pages from those that are less known. There are many derivations of PageRank \cite{langville2006updating,higham2005google}, all of which stem from modeling `websurfing' (i.e., how people navigate the web) as a Markov chain. In this analogy,   the fraction of  random websurfers at a particular web page is given by  the stationary distribution of a Markov chain. 

A main challenge for this formulation is that in practice, a network     connecting webpages via their hyperlinks does not usually consist of a single connected component. Instead, there are isolated webpages that cannot be navigated to, or navigated from. To address this issue, the Google founders introduced `teleportation' so that with probability $\alpha$, websurfers click a hyperlink to move between webpages, and with probability $1-\alpha$, websurfers randomly jump to webpage $i$ with probability $v_i$. Parameter $\alpha$ is called as \emph{teleportation parameter} or \emph{damping factor}.
In the original formulation, a walker would jump uniformly at random to another webpage so that the transition probability to each webpage is the same:  $v_i=v_j=1/N$, where $N$ is the number of webpages. It is also beneficial to allow the $v_i$ values to be heterogeneous to bias the random websurfing to remain near a particular set of webpages. Such dynamics is called `personalized' PageRank, and ${\bf v}$ is called the \emph{personalization vector}.

Both PageRank and personalized PageRank can be formulated as a discrete-time Markov chain in which the transition matrix is given by the Google matrix:
 \begin{align}
     {\bf G} = \alpha {\bf P} + (1-\alpha) {\bf e}{\bf v}^T,
 \end{align}
 where ${\bf e}$ is a vector of ones and each entry $P_{ij}$  gives the probability of transitioning from webpage $i$ to webpage $j$ following a hyperlink. The matrices $P$, ${\bf e}{\bf v}^T$ and $G$ are all row-stochastic transition matrices.
The stationary distribution of the   Markov chain with transition matrix ${\bf G}$ is called the PageRank vector $\pi\in\mathbb{R}^N$, which is the limit of the iterative equation
 \begin{align}\label{eq:it}
     {\bf x}(t+1)^T = {\bf x}(t)^T {\bf G} ,
 \end{align}
for which the fixed point is  the solution to the eigenvalue problem $ {\bf G}^T\bm{\pi} = \bm{\pi}$.

The PageRank algorithm has contributed toward Google's rise as a leading technology company, but it's worth noting that it has also been applied in a wide variety of domains beyond web search  \cite{gleich2015pagerank}. In any context, the practical usage of PageRank comes with various challenges. For example,  what $\alpha$ value should be used?   The PageRank vector can be integrated with a machine learning framework for web search \cite{najork2007hits},  used to help fit functions over directed graphs \cite{zhou2005learning}, and   can be used to predict   missing genes and protein functions \cite{ichinomiya2020protein,morrison2005generank}. In these various cases, $\alpha$ must be separately chosen as appropriate. As in \cite{pan2004automatic}, the author suggests choosing $\alpha = 0.15$ for correlated discovery in a multimedia database. Google historically set $\alpha = 0.85$, which often remains as the default choice in the literature.

\begin{figure}[t]
\begin{center}
    
\includegraphics[width=\linewidth]{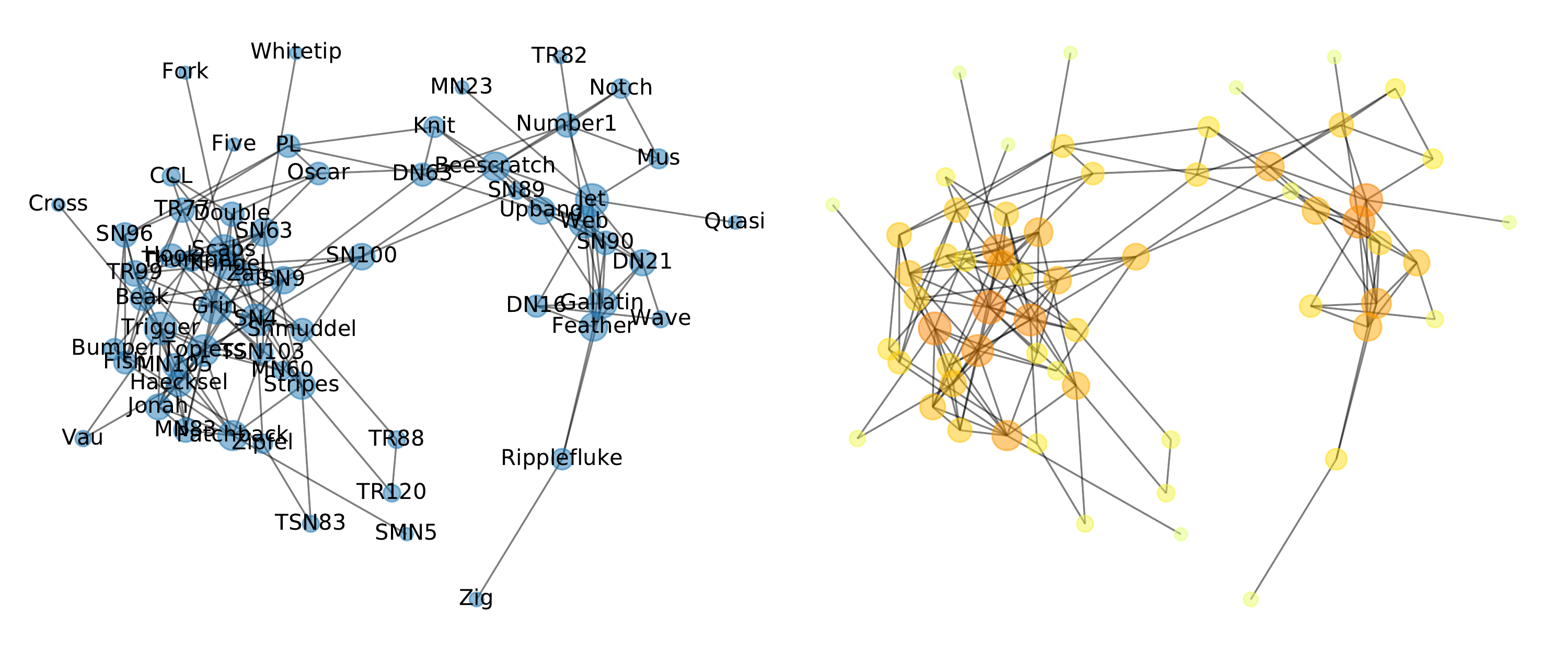}
 \caption{(A) A social network of interactions among $N=62$ dolphins in New Zealand. 
 (B) Node colors indicate the nodes' respective PageRank values. In both panels, nodes with larger/smaller size have more/fewer connecting edges.
 }
\label{Fig:Dol1}
\end{center}

\end{figure}

In our paper, we explore a different challenge that arises when considering how many times to iterate Eq.~\eqref{eq:it}. The iterated values  ${\bf x}(t)\to\bm{\pi}$ converge as $t\to\infty$, but a more practical questions involves studying how many iterations are required for the associated ranks to converge. That is, if one ranks webpages from top to bottom based on their $\pi_i$ values, then one only needs to iterate Eq.~\eqref{eq:it} until those ranks converge. Understanding the asymptotic convergence rate and asymptotic error is not immediately relevant when considering this practical question, and we propose kNN persistent homology as a mathematically principled technique to study the convergence of rankings (and relative orderings more generally).

\subsection{Ranking nodes in a dolphin social network} 
\label{sec:Dolph}
In this section, we consider the convergence of the iterative method for computing PageRank. We will study the convergence of persistence diagrams associated with converging approximate PageRank values ${\bf x}(t) \to \bm{  \pi}$, comparing the result for kNN  filtrations to those of VR filtrations. We will show that the convergence of persistence diagrams for VR filtrations closely relates to the convergence in vector norm (i.e., due to the stability theorem). In contrast,  we   the convergence of persistence diagrams for kNN filtrations more closely resembles the convergence of rank orderings. In other words, convergence of kNN persistent homology can be used to predict how many iterations are required for ${\bf x}(t)$ to be sufficiently close to $\bm{\pi}$ such that the ranking of nodes---i.e., from 1 to N---has converged.

\begin{figure}[t]
\begin{center}
    
\includegraphics[width=\linewidth]{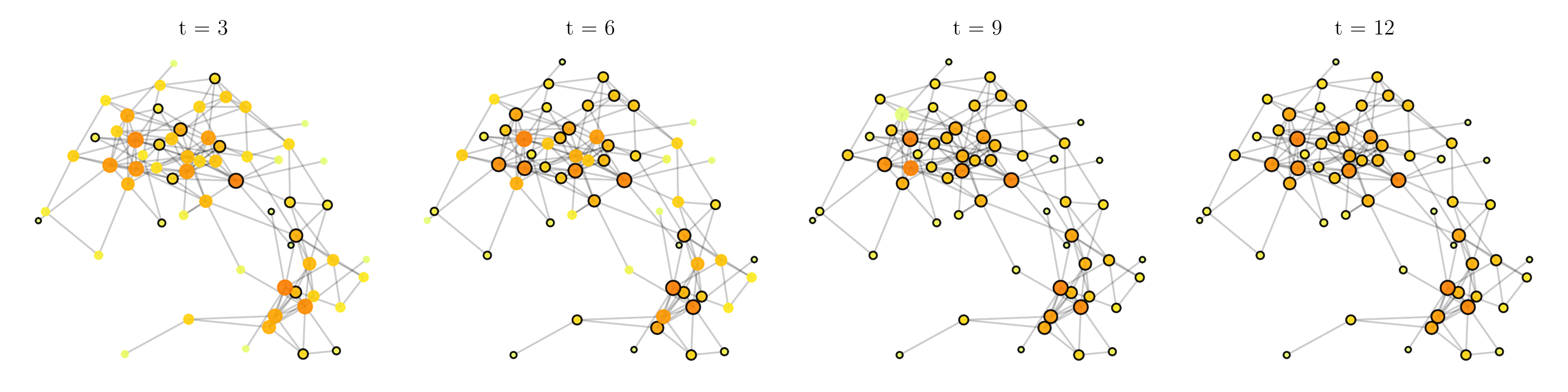}
\caption{Visualization of the approximate PageRank values $x_i(t)\approx \pi_i$ for the iterative method with different $t\in\{3,6,9,12\}$. Node colors indicate the $x_i(t)$ values. The black circles around nodes indicate the ranks  $R_i({\bf x}(t))$ that have already converged to their limiting rank  $R_i(\bm{\pi})$. For each node $i$,  we define $t^*_i$ to be the time step for which $R_i({\bf x}(t))={R}_i(\bm{\pi})$ for $t\ge t^*$. Observe that   most ranks  converge after 12 time steps, even though $||{\bf x}-\bm{\pi}|| >0$.
}
\label{Fig:Dol2}
\end{center}

\end{figure}

We present our results for an empirical network in which undirected, unweighted edges encode social interactions among $N=62$ bottlenose dolphins living near Doubtful Sound in New Zealand
\cite{lusseau2003bottlenose}. 
We apply the iterative PageRank algorithm with teleportation parameter $\alpha=0.85$ to the network with an initial condition given by $x_i(0)  = i / (\sum_{j=1}^N j)$. 
In the Figure ~\ref{Fig:Dol1}, we   illustrate the dolphin network and indicate the nodes' converged PageRank values by node color.

We formally define the converged rank order according to PageRank by
$$
[R_i(\bm{\pi)}] = \phi(-\bm{\pi}),
$$
where $\phi:\mathbb{R} \to \{1,\dots,N\}$ is the `argsort' function.
Recall that the function $\phi$ was previously defined   to sort the pairwise distances in ascending order. By multiplying $\bm{\pi}$ by negative one, we now sort the $\pi_i$ values in descending order so that $R_i(\bm{\pi})=1$ for the top-ranked node $i=\text{argmax}_j \pi_j$, which is considered to be the most important dolphin in the social network. Similarly,  $R_i(\bm{\pi})=N$ for the lowest-ranked node $i=\text{argmin}_j \pi_j$, which is   the least important dolphin. (
%
For each time $t$, we similarly define the approximate rank orderings
$$
[R_i({\bf x}(t))] = \phi(-{\bf x}(t)).
$$
We note that the rank orderings satisfy $R_i({\bf x}(t)),{R}_i(\bm{\pi}) \in \{1,\dots, N\}$.
Because ${\bf x}(t)\to\bm{\pi}$,   the approximate rank orderings $R_i({\bf x}(t))$ converge to their final rank orderings ${R}_i(\bm{\pi})$. Moreover, for each node $i$, the rank $R_i({\bf x}(t))$ can converge to $R_i(\bm{\pi})$ at a different time step $t$. Therefore, we define $t^*_i$ to be the {\bf iteration of rank convergence}, which is the value of $t$ at which $R_i({\bf x}(t))={R}_i(\bm{\pi})$ for all $t\ge t^*_i$.

In   Fig.~\ref{Fig:Dol2}, we illustrate the convergence for the approximate PageRank values $x_i(t)\approx \pi_i$ for different time steps. That is, for the time steps $t\in\{3,6,9,12\}$, we use node color to indicate the values $x_i(t)$. Moreover, we use black circles to indicate which nodes have already obtained their limiting rank, i.e., $R_i({\bf x}(t))={R(\bm{\pi})}_i$ for $t\ge t^*_i$. Most ranks converge after 12 time steps.

\begin{figure}[t]
\begin{center}
    
\includegraphics[width=\linewidth]{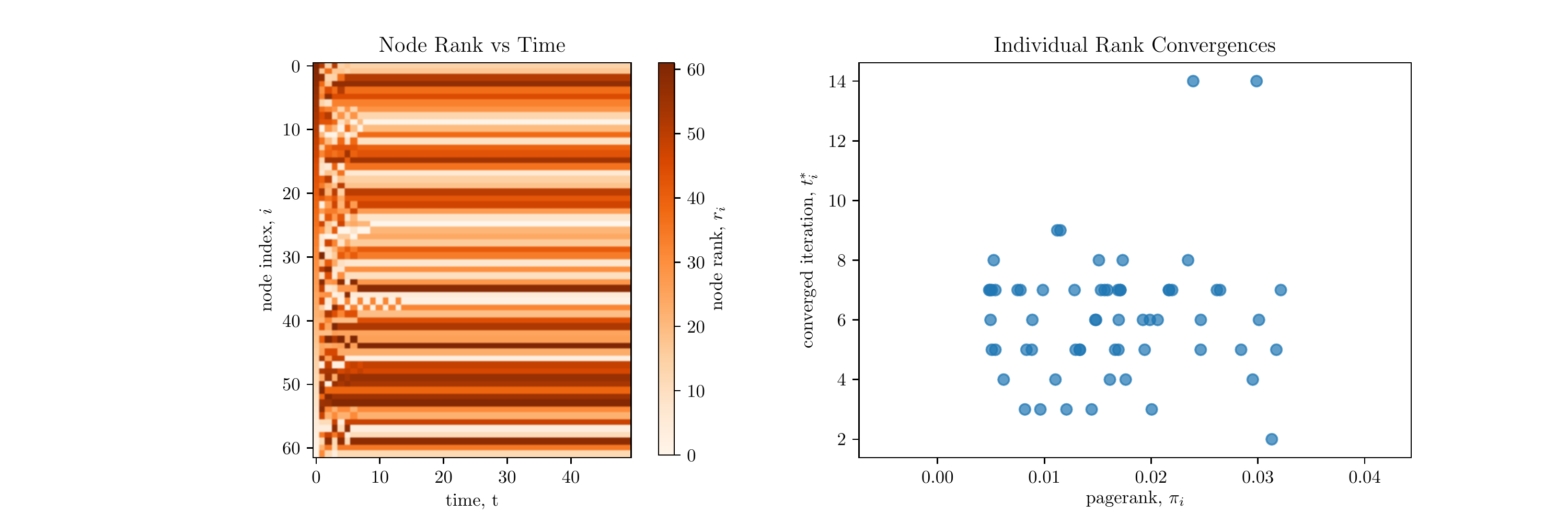}
\caption{(left) Convergence of nodes' rank orderings $R_i({\bf x}(t))\to {R}_i(\bm{\pi})$ versus time step  $t$.
(right) Scatter plot comparing $t^*_i$ and $\pi_i$ across the nodes $i$.
}
\label{Fig:BR}
\end{center}

\end{figure}

In Fig.~\ref{Fig:BR}, we further study the convergence of the rank orderings $R_i({\bf x}(t))\to {R}_i(\bm{\pi})$ for the dolphin network. In Fig.~\ref{Fig:BR}(left), observe that the $R_i({\bf x}(t))$ values change up until a time step $t^*_i$, which is potentially different  for each node $i$. In Fig.~\ref{Fig:BR}(left), we depict a scatter plot comparing $t^*_i$ and $\pi_i$, noting that we do not see any strong correlation. From a practical perspective, one is often most interested in the top-ranked nodes, and so one is primarily interested in how many iterations are required for the top-ranks to converge. However, there is no guarantee that the top-rank nodes converge before the lower-ranked ones do, or vice versa.

We now study the convergence of persistence diagrams for the  converging $x_i(t)\to\pi_i$ values, comparing persistence diagrams resulting from kNN filtrations to those of VR filtrations. More precisely, we define $D_{VR}({\bf x})$ and $D_{kNN}({\bf x})$ to be the persistence diagrams according to VR and kNN filtrations, respectively. We will also study kNN filtrations with two types of symmetrization for k-nearest neighbor sets: the minimum and maximum approaches. Given the persistence diagrams for $\bm{\pi}$ and ${\bf x}(t)$, we study homological convergence through the bottleneck distance, e.g., $d(D_{VR}({\bf x}(t))-D_{VR}(\bm{\pi}))$.

\begin{figure}[t]
\begin{center}
    
\includegraphics[width=\linewidth]{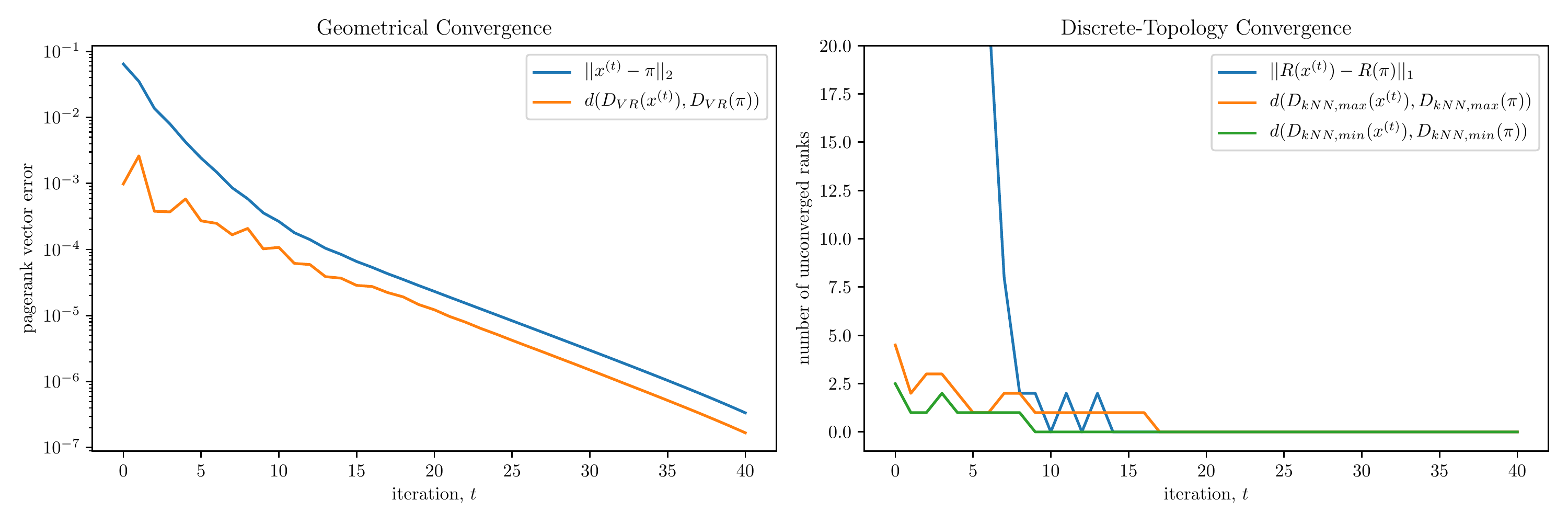}
 \caption{Homological convergence of an iterative algorithm for  PageRank for the dolphin social network.
 (left) Convergence of persistent homology for VR filtrations coincides with a geometrical notion of convergence $||{\bf x}(t)-\bm{\pi}|| \to 0$ due to the stability theorem. Both asymptotically approach 0 with exponential decay.
 (right) In contrast, convergence of persistent homology for kNN filtrations more closely resembles the convergence of the rank ordering, which exactly converges after $t=14$ time steps in this case.
 Observe the  kNN persistence diagrams for the max and min methods of symmetrization for kNN sets converge at around the same number of iterations.
  }\label{Fig:Dolconv}
  \end{center}

\end{figure}

In Fig.~\ref{Fig:Dolconv}, we illustrate the convergence of persistence diagrams for (left) VR filtrations and (right) kNN filtrations. We compare these two converging topological spaces, respectively, with the normed approximation error, $||{\bf x}(t)-\bm{\pi}||$, and the total difference in rank orderings, $\sum_i |{R}_i(\bm{\pi}) -  R_i({\bf x}(t))|$. Observe in Fig.~\ref{Fig:Dolconv}(left) that VR persistent homology converges similarly to  $||{\bf x}(t)-\bm{\pi}||$, whereas kNN persistent homology converges similarly to $\sum_i |{R}_i(\bm{\pi}) -  R_i({\bf x}(t))|$. That is, kNN topological convergence can be used as a proxy to estimate the number of iterations required for the   rank orderings to converge to exactly their final values.

\subsection{U-local convergence of dolphin social network} \label{sec:Dolpu}

Before concluding, we further study the homological convergence of PageRank for the dolphin network, except we now focus on 
our notion of $\mathcal{U}$-local topological convergence, which we presented in Sec.~\ref{sec:convergence1}. In this context, we focus on convergence for a subset $\mathcal{U}\subset \mathcal{V}$ of the nodes, and we will consider two subsets: a randomly selected set of nodes---$\mathcal{U}_1= \{$Ripplefluke, Zig, Feather, Gallatin, SN90, DN16, Wave, DN21, Web, Upbang$\}$ and $\mathcal{U}_2$ is the set of 10  nodes that have the  top PageRank values. Note that we still compute the iterative approximation to PageRank in the usual way, except that we only consider the values $\pi_i$ and $x_i(t)$ for which $i\in\mathcal{U}$.

\begin{figure}[t]
\begin{center}
    
\includegraphics[width=\linewidth]{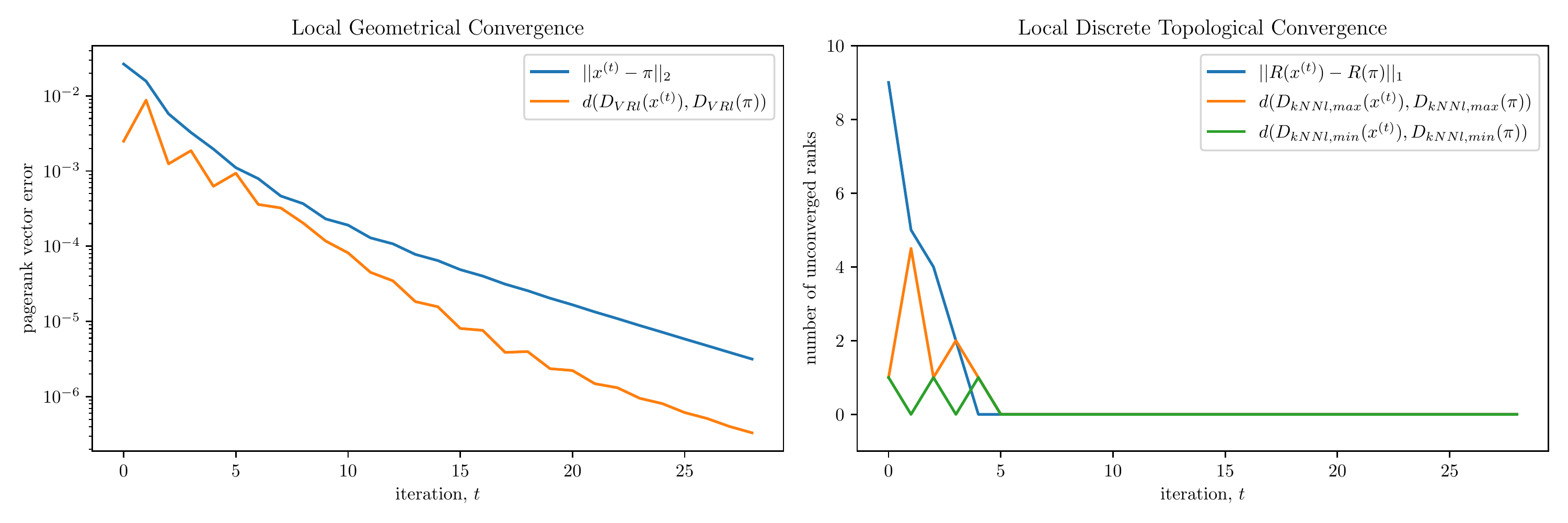}
 \caption{$\mathcal{U}$-local topological convergence for the following subset of nodes: $\mathcal{U}_1=\{$Ripplefluke, Zig, Feather, Gallatin, SN90, DN16, Wave, DN21, Web, Upbang$\}$ . Similar to Fig.~\ref{Fig:Dolconv}, the left and right panels depict convergence of persistence diagrams for VR and kNN filtrations, respectively.
  }\label{Fig:locd}
  \end{center}

\end{figure}

\begin{figure}[t]
\begin{center}
    
\includegraphics[width=\linewidth]{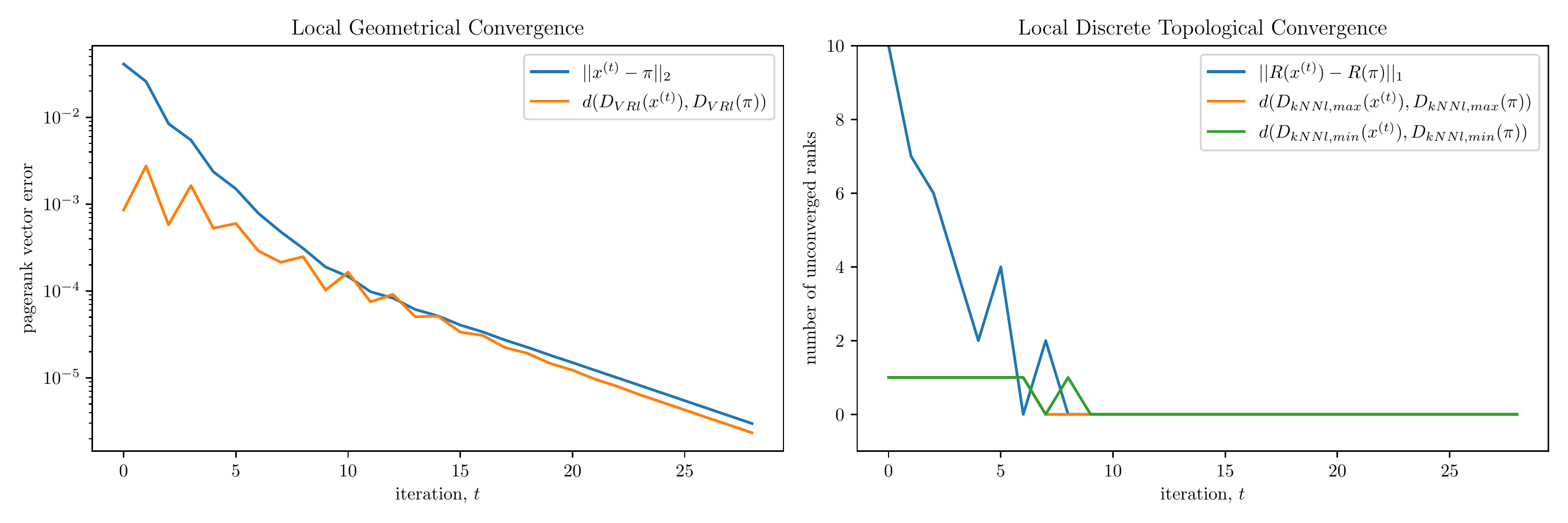}
 \caption{Same information as in Fig.~\ref{Fig:locd}, except we consider  the subset $\mathcal{U}_1$ of nodes that have the largest PageRank values.
 }\label{Fig:top10d}
 \end{center}

\end{figure}

In Fig.~\ref{Fig:locd}(left) and (right), we depict convergence of VR and kNN persistence homology, respectively for the subset of nodes $\mathcal{U}_1$.  Observer that convergence for VR persistent homology is similar to that for the full system. In contrast, observe in Fig.~\ref{Fig:locd}(right) that the kNN homology and rank orderings converge after approximately $t=5$ iterations for this subset.

In Fig.~\ref{Fig:top10d}, we depict the same information, except that we now consider the subset $\mathcal{U}_2$ of nodes with top PageRank values. 
In this case, the VR convergence more closely aligns with $||{\bf x}(t)-\bm{\pi}||$ than when we considered the full set $\mathcal{V}$ of nodes or the subset $\mathcal{U}_1$. Moreover, the kNN homology and the rank orderings converge after approximately $t=9$ iterations in this case.

\section{Discussion}\label{sec:discusion}

In this paper, we developed an approach for topological data analysis (TDA) that utilizes k-nearest neighbor sets to define kNN complexes, kNN filtrations, and kNN persistent homology. Our approach was developed 
in the spirit of discrete topology and by examining the relative ordering of points, as opposed to the precise distances between points. 
Although kNN orderings are related to pairwise distances, we provided theory and many experiments highlighting important differences between  persistent homology that is based on kNN filtrations versus Vietoris-Rips (VR) filtrations.

To gain theoretical insights into kNN-based TDA, 
we investigated   stability properties for the resulting persistence diagrams
in Sec.~\ref{sec: stability1} and convergence properties in
Sec.~\ref{sec:convergence1}.
While persistence diagrams resulting from kNN filtrations do not satisfy a  stability theorem involving a universal bound on perturbed point sets (see Theorem~\ref{thm:kNN_sta}), we identified and characterized different notions of stability and convergence by identifying equivalence classes as well as bounds on the bottleneck distance between persistence diagrams for kNN homology that are based on the maximum difference for a kNN ordering. Our formulation of convergence also led to several types  including global convergence, $\kappa$-bounded convergence and $\mathcal{U}$-local convergence. Moreover, we identified some relations among these types as in Theorem~\ref{thm: convergence_NNO} and Theorem~\ref{thm: convergence NNO1}.

As a concrete application, we applied persistent homology to study the topological convergence of an iterative method for solving PageRank. We showed that  the convergence of persistence diagrams for VR filtrations closely relates to the convergence in vector norm (i.e., due to the stability theorem). In contrast,    the convergence of persistence diagrams for kNN filtrations more closely resembles the convergence of rank ordering. In other words, convergence of kNN persistent homology can be used to predict how many iterations are required for ${\bf x}(t)$ to be sufficiently close to $\bm{\pi}$ such that the ranking of nodes---i.e., from 1 to N---has converged. Although we have focused on the PageRank algorithm, iterative algorithms for solving systems (e.g., root finding) are some of the most widely used numerical algorithms. We have shown that the existing TDA approach of VR filtrations coincides with geometrical (i.e., normed) convergence and thereby provides one perspective for the convergence of numerical algorithms. In contrast, kNN filtrations provide complementary insights from the perspective of a discrete topological space that is associated with the relative positioning of points. As such, we expect kNN complexes, filtrations, and persistent homology to have many applications for converging  numerical algorithms beyond PageRank.

\section*{Acknowledgments} M.Q.L. and D.T. were supported in part by the Simons Foundation (Grant No. 578333). D.T also acknowledges  the National Science Foundation (Grants No. DMS-2052720 and No. EDT-1551069).








\bibliographystyle{plain}
\bibliography{ref}

\medskip
Received xxxx 20xx; revised xxxx 20xx; early access xxxx 20xx.
\medskip

\end{document}